\documentclass[10 pt]{amsart}
\usepackage{amsmath}
\usepackage{amssymb}
\usepackage{graphicx}
\usepackage{amsthm}
\usepackage{stmaryrd}
\usepackage[all]{xy}

\newtheorem{theorem}{Theorem}
\newtheorem{proposition}[theorem]{Proposition}
\newtheorem{lemma}[theorem]{Lemma}

\newtheorem{claim}{Claim}
\theoremstyle{definition}
\newtheorem{definition}{Definition}

\newtheorem{fact}{Fact}

\makeatletter

\def\dotminussym#1#2{%
  \setbox0=\hbox{$\m@th#1-$}%
  \kern.5\wd0%
  \hbox to 0pt{\hss\hbox{$\m@th#1-$}\hss}%
  \raise.35\ht0\hbox to 0pt{\hss$\m@th#1\cdot$\hss}
  \kern.5\wd0}
\newcommand{\dotminus}{\mathbin{\mathpalette\dotminussym{}}}

\makeatother

\title{Vaught's Two-Cardinal Theorem and Quasi-Minimality in Continuous Logic}
\date{\today}

\author{Victoria Noquez}
\address{Indiana University Bloomington, Department of Mathematics, 831 East 3rd St., Bloomington, IN 47405-7106}
\email{vnoquez@iu.edu}

\begin{document}

\maketitle 

\begin{abstract} 
We prove the following continuous analogue of Vaught's Two-Cardinal Theorem: if for some $\kappa>\lambda\geq \aleph_0$, a continuous theory $T$ has a model with density character $\kappa$ which has a definable subset of density character $\lambda$, then $T$ has a model with density character $\aleph_1$ which has a separable definable subset.  We also show that if we assume that $T$ is $\omega$-stable, then if $T$ has a model of density character $\aleph_1$ with a separable definable set, then for any uncountable $\kappa$ we can find a model of $T$ with density character $\kappa$ which has a separable definable subset.  In order to prove this, we develop an approximate notion of quasi-minimality for the continuous setting.  We apply these results to show a continuous version of the forward direction of the Baldwin-Lachlan characterization of uncountable categoricity: if a continuous theory $T$ is uncountably categorical, then $T$ is $\omega$-stable and has no Vaughtian pairs.  
\end{abstract}

\section{Introduction}

In classical logic, the study of uncountable categoricity has led to the development of numerous model theoretic tools, many of which have been very useful for a variety of different applications.  

There has been significant progress towards understanding uncountable categoricity in the continuous setting, particularly for examples from functional analysis. In \cite{uands}, Usvyatsov and Shelah prove that every model of an uncountably categorical theory expanding a Banach space is prime over a spreading model, isometric to the standard basis of a Hilbert space.  More recently, in \cite{henson}, Henson and Raynaud provide a criterion for ensuring that the elementary class of a modular Banach space consists of all direct sums of that space with arbitrary Hilbert spaces.  This leads to many examples of uncountably categorical theories whose models are Banach spaces. 

Towards understanding uncountable categoricity in continuous logic outside of the context of Banach spaces, in \cite{benyaacov}, Ben Yaacov proves Morley's categoricity theorem for continuous logic: for $\kappa,\lambda>\aleph_0$, a continuous theory $T$ is $\kappa$-categorical if and only if $T$ is $\lambda$-categorical.  

In the classical setting, Morley's categoricity theorem can be viewed as a corollary of the following theorem of Baldwin and Lachlan: for $\kappa$ an uncountable cardinal, a theory $T$ is $\kappa$-categorical if and only if $T$ is $\omega$-stable and has no Vaughtian pairs.  

The work in this paper was motivated by an effort to show a continuous version of the Baldwin-Lachlan characterization of uncountable categoricity.  We resolve the forward direction by showing Vaught's Two-Cardinal Theorem as well as a partial converse to the theorem which requires an assumption of $\omega$-stability.  

We begin by discussing definability in continuous logic.  Recall that in the continuous setting, formulas in the language are functions to $[0,1]$.  Definable sets are not nearly as well behaved as they are in the classical setting, and it is important in this context that we carefully distinguish between zero sets of definable predicates (functions which can be uniformly approximated by formulas) and definable sets themselves.  We will call a set $D$ \emph{definable} if the predicate $dist(x,D) = \displaystyle{\inf_{y\in D}} d(x,y)$ can be uniformly approximated by formulas in the language.  Note that zero sets of definable predicates may not be definable sets.  Throughout this paper, we adopt the convention of saying that a definable predicate $P$ is a \emph{definable distance predicate} if and only if $Z(P)$ is a definable set.

In Section 3 we give a continuous definition of a Vaughtian pair: $(\mathcal{N},\mathcal{M})$ is a Vaughtian pair of models of a continuous theory $T$ if $\mathcal{M}\precneqq \mathcal{N}$ and there is a definable distance predicate $P$ (the zero set of $P$ is a definable set) such that $\{x\in\mathcal{M}|P(x) = 0\} = \{x\in\mathcal{N}|P(x)=0\}$. We describe how Vaughtian pairs can be viewed as structures in an extension of our language and lay the groundwork necessary to show that if $T$ has a Vaughtian pair of models, then $T$ has a Vaughtian pair of models $(\mathcal{N},\mathcal{M})$ such that $\mathcal{M}$ and $\mathcal{N}$ are separable and $\mathcal{M}\cong\mathcal{N}$. The work in this section is mostly analogous to the work on Vaughtian pairs in the classical setting, with minor adjustments being made to adapt to the continuous setting.

In Section 4, we introduce a continuous version of quasi-minimality and see that it requires a significantly different approach than in the classical setting.  This is by far the most technical section.  The purpose of defining quasi-minimality in the continuous setting is to prove Lemma \ref{qmcor}, which in turn allows us to prove Theorem \ref{kappamodel} (described in more detail shortly), which we need for our main result about uncountable categoricity.   This section is where we see the greatest departure from the classical setting in terms of the mechanics of the proof.

In Section 5, we introduce $(\kappa,\lambda)$-models: for a continuous theory $T$, and $\kappa>\lambda\geq \aleph_0$, $T$ has a $(\kappa,\lambda)$-model if there is $\mathcal{M}\vDash T$ such that $\mathcal{M}$ has density character $\kappa$ and a definable subset with density character $\lambda$.  The existence of such a model is useful in understanding uncountable categoricity: by a compactness argument, we can show that for uncountable $\kappa$, a theory $T$ will have a model of density character $\kappa$ such that every non-compact definable set has density character $\kappa$.  So if $T$ has a $(\kappa,\lambda)$-model with $\kappa>\lambda\geq \aleph_0$, it cannot be $\kappa$-categorical.

Theorem \ref{vaughtstwocardinaltheorem} is the continuous analogue of Vaught's Two-Cardinal theorem: for $\kappa>\lambda\geq \aleph_0$, if a continuous theory $T$ has a $(\kappa,\lambda)$-model, then $T$ has an $(\aleph_1,\aleph_0)$-model.   We prove this in an analogous way to its classical counterpart.  

In Theorem \ref{kappamodel}, we  show that if $T$ is assumed to be $\omega$-stable, then if $T$ has an $(\aleph_1,\aleph_0)$-model, then for any uncountable $\kappa$, $T$ has a $(\kappa,\aleph_0)$-model.   Using the work from Section 4, we are able to proceed as in the classical case (though with some extra technical considerations).  

In Section 6, we apply the results of Section 5 to show that if a continuous theory $T$ is $\kappa$-categorical for some uncountable $\kappa$, then it has no Vaughtian pairs.  This along with Theorem 5.2 in \cite{benyaacov} (if $T$ is uncountably categorical, it is $\omega$-stable) gives us the forward direction of the Baldwin-Lachlan characterization of uncountable categoricity in continuous logic. Towards understanding the converse, it was shown in \cite{noquez} that an argument analogous to the classical argument will fail, because for a reasonable notion of minimality, having no Vaughtian pairs of models of a theory does not guarantee that a minimal definable predicate is strongly minimal.  This question was then fully resolved in \cite{hanson2020strongly}, in which Hanson gave an example of an $\omega$-stable theory with no Vaughtian pairs which fails to be inseparably categorical. 

In Section 7, we provide two examples of Vaughtian pairs of models.  The first is a Vaughtian pair of models of the theory of the Urysohn sphere, which we obtain by removing a small portion of the sphere, and then describing a definable set which is sufficiently far away from the removed portion.  The second example is a Vaughtian pair of models of the randomization of any (continuous or classical) theory $T$ is a countable language.  Though we can see this directly via product randomizations (as in \cite{AK}), by Theorem \ref{forwarddirection}, this shows us that uncountable categoricity is not preserved by randomizations.

\section{Definability in Continuous Logic}

We assume that the reader is familiar with the basics of continuous model theory.  An in depth introduction to the topic can be found in \cite{mtfms}.   For the reader's convenience, we will discuss definability in the continuous setting. 


Let $\mathcal{L}$ be a continuous language and $\mathcal{M}$ an $\mathcal{L}$-structure. 
A function $P:\mathcal{M}^n\rightarrow[0,1]$ is a \emph{definable predicate} in $\mathcal{M}$ over some $A\subset\mathcal{M}$ if there exists a sequence $(\phi_k(x):k<\omega)$ of $\mathcal{L}(A)$-formulas such that $\phi_k^{\mathcal{M}}(x)$ converge uniformly to $P(x)$ on all of $\mathcal{M}^n$. In this case we say that the predicate $P$ is definable over $A$.  The \emph{zero set} of $P$ in $\mathcal{M}$, $Z(P^{\mathcal{M}}) = \{x\in\mathcal{M}|P(x)=0\}$.   If $\mathcal{M}$ is clear from context, we just write $Z(P)$.   



For $D\subset\mathcal{M}^n$, we say that the set $D$ is \emph{definable} if $dist(x,D) = \displaystyle{\inf_{y\in D}}\ \displaystyle{\max_{1\leq i\leq n}} d(x_i,y_i)$ is a definable predicate.  So $D = Z(dist(x,D))$.  

If $P$ is a definable predicate such that $Z(P)$ is a definable set, we call $P$ a  \emph{definable distance predicate}.  Thus, our definable sets are exactly the zero sets of definable distance predicates.

We have the following equivalent conditions for definability of sets (Theorem 9.17 and Proposition 9.19 in \cite{mtfms}):

\begin{proposition}\label{defdistpred}  Let $D$ be a set in some $\mathcal{L}$-structure, $\mathcal{M}$.  The following are equivalent:
\begin{enumerate}
\item $D$ is definable in $\mathcal{M}$ over $A$.
\item For any definable predicate $Q(x,y)$, $\displaystyle{\sup_{y\in D}} Q(x,y)$ and $\displaystyle{\inf_{y\in D}} Q(x,y)$ are definable in $\mathcal{M}$ over $A$.
\item There is a predicate $P:\mathcal{M}^n\rightarrow [0,1]$ definable in $\mathcal{M}$ over $A$, such
that $D\subset Z(P)$ and $\forall \epsilon>0\exists \delta>0 \forall x\in\mathcal{M}(P(x)\leq \delta\Rightarrow dist(x,D)\leq \epsilon)$.  So $D=Z(P)$.  
\end{enumerate}
\end{proposition}

\section{Vaughtian Pairs}

In order to prove Vaught's Two-Cardinal Theorem, we must develop a continuous analogue of Vaughtian pairs.  We propose the following definition of Vaughtian pairs.

 \begin{definition} Let $T$ be a continuous $\mathcal{L}$-theory.  $(\mathcal{N},\mathcal{M})$ is a \emph{Vaughtian pair} of models of $T$ if $\mathcal{M}\prec\mathcal{N}$, $\mathcal{M}\neq\mathcal{N}$, and there is $P$, a definable distance predicate over $\mathcal{M}$, such that $Z(P^{\mathcal{M}})$ is not compact, and $Z(P^{\mathcal{M}}) = Z(P^{\mathcal{N}})$.  
\end{definition}

Note that in the continuous setting, compact is analogous to finite.  

In defense of this particular characterization, as opposed to, for example, just considering zero sets of definable predicates which are not necessarily definable sets themselves, in Section 7 we provided examples of Vaughtian pairs of models in the continuous setting.  

The proofs in the remainder of this section are analogous to those in the classical setting, with only small adjustments necessary to adapt to continuous logic.

Let $\mathcal{L}^* = \mathcal{L}\cup\{U\}$ where $U$ is a unary predicate.  We can view a pair $(\mathcal{N},\mathcal{M})$ of models of $T$ as an $\mathcal{L}^*$-structure by interpreting $U(x) = C\cdot \displaystyle{\inf_{y\in\mathcal{M}}} d(x,y)$, where $C>0$ is some constant (so $Z(U) = \mathcal{M}$, and $U$ is a definable distance predicate, since $d(x,Z(U)) = \frac{1}{C}U(x)$).  

For an $\mathcal{L}$-formula $\phi(\overline{x})$, we will define $\phi^U(\overline{x})$, the restriction of $\phi$ to $Z(U)$, inductively as follows:
\begin{itemize}
\item If $\phi(\overline{x})$ is atomic, then $\phi^U(\overline{x}) = \max(U(x_1),\ldots,U(x_n),\phi(\overline{x}))$.
\item $(\frac{1}{2}\phi(\overline{x}))^U = \frac{1}{2} \phi^U(\overline{x})$.
\item $(\phi(\overline{x})\dotminus \psi(\overline{x}))^U = \phi^U(\overline{x})\dotminus \psi^U(\overline{x})$.  
\item $(\displaystyle{\sup_x}\ \phi(x,\overline{a}))^U = \displaystyle{\sup_{x\in Z(U)}} \phi^U(x,\overline{a})$.  
\item $(\displaystyle{\inf_x}\ \phi(x,\overline{a}))^U = \displaystyle{\inf_{x\in Z(U)}} \phi^U(x,\overline{a})$. 
\end{itemize}

Observe that by Proposition \ref{defdistpred}, $(\displaystyle{\sup_x} \phi(x,\overline{a}))^U$ and $(\displaystyle{\inf_x} \phi(x,\overline{a}))^U$ are definable predicates for any $\mathcal{L}$-formula $\phi$.

An easy induction shows that for each formula $\phi$, for all $r\in[0,1]$ and $\overline{a}\in \mathcal{M}$, 

$\phi^U(\overline{a})\leq r\Leftrightarrow \mathcal{M}\vDash\phi(\overline{a})\leq r$, and $\phi^U(\overline{a})\geq r\Leftrightarrow \mathcal{M}\vDash\phi(\overline{a})\geq r$.  

Thus, $\phi^U(\overline{a})= r\Leftrightarrow \mathcal{M}\vDash\phi(\overline{a})= r$.

\begin{lemma}\label{countablevp} If $T$ has a Vaughtian pair $(\mathcal{N},\mathcal{M})$, then there is a Vaughtian pair $(\mathcal{N}_0,\mathcal{M}_0)$ where $\mathcal{N}_0$ has countable density character.
\begin{proof}
Let $(\mathcal{N},\mathcal{M})$ be a Vaughtian pair of models of $T$.  Let $P$ be a definable distance predicate over some countable $A\subset\mathcal{M}$ such that $Z(P^{\mathcal{M}})$ is not compact and $Z(P^{\mathcal{M}}) = Z(P^{\mathcal{N}})$.   Since $(\mathcal{N},\mathcal{M})$ is a Vaughtian pair, there is $x\in\mathcal{N}\setminus\mathcal{M}$ with $d(x,\mathcal{M})=\delta>0$, so view $(\mathcal{N},\mathcal{M})$ as an $\mathcal{L}^*$-structure by interpreting $U(z)$ as $\max(1,\frac{1}{\delta}\cdot d(z,\mathcal{M}))$.  

By downward L{\"o}wenheim-Skolem (Proposition 7.3 in \cite{mtfms}), there is $(\mathcal{N}_0,\mathcal{M}_0)$ with $A\subset\mathcal{N}_0$ where $\mathcal{N}_0$ has countable density character and $(\mathcal{N}_0,\mathcal{M}_0)\prec (\mathcal{N},\mathcal{M})$ as $\mathcal{L}^*_A$-structures.  

$(\mathcal{N},\mathcal{M})$ satisfies $U(c)=0$ for any constants $c$ in $\mathcal{L}$ and $\displaystyle{\sup_{\overline{x}\in Z(U)}} U(f(\overline{x})) = 0$ for function symbols $f$, where $\overline{x}\in Z(U)$ is shorthand for 

$\overline{x}\in Z(\max(U(x_1),\ldots,U(x_n)))$.  Thus, $(\mathcal{N}_0,\mathcal{M}_0)$ satisfies these as well, so 

\noindent $\mathcal{M}_0 = Z(U^{(\mathcal{N}_0,\mathcal{M}_0)})$ is a substructure of $\mathcal{N}_0$.  

Since $\mathcal{M}\prec\mathcal{N}$, $(\mathcal{N},\mathcal{M})\vDash\displaystyle{\sup_{\overline{x}\in Z(U)}} |\phi^U(\overline{x}) - \phi(\overline{x})| = 0$ for all $\mathcal{L}$-formulas $\phi(\overline{x})$, so $(\mathcal{N}_0,\mathcal{M}_0)\vDash \displaystyle{\sup_{\overline{x}\in Z(U)}} |\phi^U(\overline{x}) - \phi(\overline{x})| = 0$, which means that $\mathcal{M}_0\prec\mathcal{N}_0$.  Also note that since $(\mathcal{N},\mathcal{M})\vDash U(a)=0$ for all $a\in A$, $(\mathcal{N}_0,\mathcal{M}_0)\vDash U(a)=0$ for all $a\in A$, so $A\subset\mathcal{M}_0$.  

To see that $\mathcal{N}_0\neq\mathcal{M}_0$, we know that $(\mathcal{N},\mathcal{M})\vDash \displaystyle{\sup_{x}} |U(x)-1| = 0$, so 

\noindent $(\mathcal{N}_0,\mathcal{M}_0)\vDash \displaystyle{\sup_{x}} |U(x)-1| = 0$, which means there exists $x\in \mathcal{N}_0$ such that $U(x)\neq 0$ (since there exists $x\in\mathcal{N}_0$ such that $U(x)$ is arbitrarily close to $1$).  So $x\in\mathcal{N}_0\setminus\mathcal{M}_0$.  

Since $P$ is a definable distance predicate and 

\noindent $(\mathcal{N},\mathcal{M})\vDash\displaystyle{\sup_{\overline{x}\in Z(P)}} \max(U(x_1),\ldots,U(x_n)) = 0$, we must have  

\noindent $(\mathcal{N}_0,\mathcal{M}_0)\vDash\displaystyle{\sup_{\overline{x}\in Z(P)}} \max(U(x_1),\ldots,U(x_n)) = 0$ which means that $Z(P^{\mathcal{N}_0}) = Z(P^{\mathcal{M}_0})$.  

Since $Z(P^{\mathcal{N}})$ is not totally bounded, for some $\epsilon>0$ there is no finite $\epsilon$-net.  Suppose for some $n<\omega$, $x_1,\ldots,x_{n}\in\mathcal{N}_0$ are the centers of a finite $\epsilon$-net of $Z(P^{\mathcal{N}_0})$.  Then $\mathcal{N}_0\vDash \displaystyle{\sup_{y\in Z(P)}}\ \displaystyle{\min_{1\leq i\leq n}}\ d(y,x_i) \dotminus \epsilon=0$, so by elementarily, 

\noindent $\mathcal{N}\vDash  \displaystyle{\sup_{y\in Z(P)}}\ \displaystyle{\min_{1\leq i\leq n}}\ d(y,x_i) \dotminus \epsilon=0$, so this is a finite $\epsilon$-net of $Z(P^{\mathcal{N}})$.  Thus, $\mathcal{N}_0$ cannot have a finite $\epsilon$-net.  So $Z(P^{\mathcal{N}_0})$ is not totally bounded.


Hence, $(\mathcal{N}_0,\mathcal{M}_0)$ is a Vaughtian pair.  
\end{proof}
\end{lemma}

\begin{proposition}\label{sametypes} Let $T$ be a complete continuous theory in a countable language $\mathcal{L}$.  If $\mathcal{M}$ and $\mathcal{N}$ are homogeneous models of $T$ with countable density character which realize the same types in $S_n(T)$ for $n\geq 1$, then $\mathcal{M}\cong \mathcal{N}$. 
\begin{proof}
Let $\mathcal{M}_0\subset\mathcal{M}$ and $\mathcal{N}_0\subset\mathcal{N}$ be countable dense subsets of $\mathcal{M}$ and $\mathcal{N}$ respectively.  We will build a bijective elementary map $f:\widehat{\mathcal{M}_0}\rightarrow\widehat{\mathcal{N}_0}$ where $\mathcal{M}_0\subset\widehat{\mathcal{M}_0}\subset\mathcal{M}$ and $\mathcal{N}_0\subset\widehat{\mathcal{N}_0}\subset\mathcal{N}$ by a back-and-forth argument.  Then we can uniquely extend this to an isomorphism from $\mathcal{M}\rightarrow\mathcal{N}$: for $x\in\mathcal{M}\setminus\widehat{\mathcal{M}_0}$, let $x_i\in\widehat{\mathcal{M}_0}$ such that $x_i\rightarrow x$.  Then, since $(x_i:i<\omega)$ is Cauchy and $f$ is elementary, $(f(x_i):i<\omega)$ is Cauchy, and thus, converges to some $y\in\mathcal{N}$.  If $y\in \widehat{\mathcal{N}_0}$, there is $x'\in\widehat{\mathcal{M}_0}$ such that $f(x')=y$.  So since $\forall \epsilon>0$ there is $N<\omega$ such that for all $i\geq N$, $\mathcal{N}\vDash d(f(x_i),f(x'))\dotminus  \epsilon =0$, by the elementarity of $f$, $\mathcal{M}\vDash d(x_i,x')\dotminus \epsilon=0$.  So $x_i\rightarrow x'$, which means $x=x'$, contradicting that $x\notin\widehat{\mathcal{M}_0}$.  This $y$ is unique, as it is the limit of $f(x_i)$, and $f\cup \{(x,y)\}$ is elementary, since $f$ is elementary, and thus, continuous.  


We will build a sequence of partial elementary maps $f_0\subset f_1\subset\ldots$ with finite domain and let $f = \displaystyle{\bigcup_{i<\omega}} f_i$, so $f$ will be elementary.

Let $a_i$ enumerate $\mathcal{M}_0$ and $b_i$ enumerate $\mathcal{N}_0$.  We will ensure that $a_i\in dom(f_{2i+1})$ and $b_i\in img(f_{2i+2})$.  Thus, we will have $\mathcal{M}_0\subset dom(f)=\widehat{\mathcal{M}}_0$ and $f:\widehat{\mathcal{M}_0}\rightarrow\widehat{\mathcal{N}_0}$ surjective where $\widehat{\mathcal{N}_0}\supset \mathcal{N}_0$.  

Stage 0: $f_0=\emptyset$.  This is partial elementary since $T$ is complete.

Inductively assume that $f_s$ is partial elementary.  Let $\overline{a}$ be the domain of $f_s$ and $\overline{b} = f_s(\overline{a})$.  

Stage $s+1=2i+1$:  If $a_i$ is in $\overline{a}$, $f_{s+1} = f_s$.  Suppose not.  Since $\mathcal{M}$ and $\mathcal{N}$ realize the same types, there is $\overline{c},d\in\mathcal{N}$ such that $tp^{\mathcal{N}}(\overline{c},d) = tp^{\mathcal{M}}(\overline{a},a_i)$.  So in particular, $tp^{\mathcal{N}}(\overline{c}) = tp^{\mathcal{M}}(\overline{a})=tp^{\mathcal{N}}(\overline{b})$ since $f_s$ is elementary.  Since $\mathcal{N}$ is homogeneous, there is $e\in\mathcal{N}$ such that $tp^{\mathcal{N}}(\overline{b},e) = tp^{\mathcal{N}}(\overline{c},d) = tp^{\mathcal{M}}(\overline{a},a_i)$.  So $f_{s+1} = f_s\cup\{(a_i,e)\}$ is partial elementary and has $a_i$ in its domain.  

Stage $s+1=2i+2$: Again, if $b_i\in img(f_s)$, then $f_{s+1} =f _s$.  If not, as in the previous stage, there are $\overline{c},d\in\mathcal{M}$ such that $tp^{\mathcal{M}}(\overline{c},d) = tp^{\mathcal{N}}(\overline{b},b_i)$, so $tp^{\mathcal{M}}(\overline{c}) = tp^{\mathcal{N}}(\overline{b}) = tp^{\mathcal{M}}(\overline{a})$ since $f_s$ is elementary.  Since $\mathcal{M}$ is homogeneous, there is $e\in\mathcal{M}$ such that $tp^{\mathcal{M}}(\overline{a},e) = tp^{\mathcal{M}}(\overline{c},d) = tp^{\mathcal{N}}(\overline{b},b_i)$.  So $f_{s+1} = f_s\cup\{(e,b_i)\}$ is elementary and has $b_i$ in its image.  
\end{proof}
\end{proposition}

\begin{lemma}\label{homogvp} Suppose $\mathcal{M}_0\precneqq\mathcal{N}_0$ are models of $T$ with countable density character.  We can find $\mathcal{M},\mathcal{N}\vDash T$ such that  $(\mathcal{N}_0,\mathcal{M}_0)\prec(\mathcal{N},\mathcal{M})$ (as $\mathcal{L}^*$-structures) and so that $\mathcal{N}$ and $\mathcal{M}$ have countable density character, are homogeneous, and realize the same types in $S_n(T)$ for $n\geq 1$.  

So by the previous proposition, $\mathcal{M}\cong\mathcal{N}$.
\begin{proof} First we need a few claims.  

\begin{claim} If $\overline{a}\in\mathcal{M}_0$ and $p\in S_n(\overline{a})$ is realized in $\mathcal{N}_0$, then there is $(\mathcal{N}',\mathcal{M}')\succ (\mathcal{N}_0,\mathcal{M}_0)$ with countable density character such that $p$ is realized in $\mathcal{M}'$.  
\begin{proof}
Consider the type $\{\phi^U(\overline{x},\overline{a})\dotminus \frac{1}{n}|\phi(\overline{x},\overline{a})\in p,n<\omega\}\cup Diag_{el}(\mathcal{N}_0,\mathcal{M}_0)$.  Let $\Gamma(\overline{x})$ be a finite subset of this and let $\phi_1(\overline{x},\overline{a}),\ldots,\phi_m(\overline{x},\overline{a})\in p$ be such that 

$\Gamma(\overline{x})\subset\{\phi_1^U(\overline{x},\overline{a})\dotminus \frac{1}{n_1},\ldots,\phi_m^U(\overline{x},\overline{a})\dotminus \frac{1}{n_m}\}\cup Diag_{el}(\mathcal{N}_0,\mathcal{M}_0)$ for some 

$n_1,\ldots,n_m<\omega$.  $\mathcal{N}_0\vDash \displaystyle{\inf_{\overline{x}}} \displaystyle{\max_{1\leq i\leq m}} \phi_i(\overline{x},\overline{a})=0$ since $p$ is realized in $\mathcal{N}_0$, so $\mathcal{M}_0\vDash \displaystyle{\inf_{\overline{x}}} \displaystyle{\max_{1\leq i\leq m}} \phi_i(\overline{x},\overline{a}) =0$.  Let $n=\displaystyle{\max_{1\leq i\leq m}} n_i$ and choose $\overline{x}\in\mathcal{M}_0$ such that $\displaystyle{\max_{1\leq i\leq m}} \phi_i(\overline{x},\overline{a}) \leq \frac{1}{n}$.  $\overline{x}$ will realize $\Gamma(\overline{x})$.  So since the type is finitely satisfiable, there is $(\mathcal{N}',\mathcal{M}')$ with countable density character realizing it.  Since it includes the elementary diagram of $(\mathcal{N}_0,\mathcal{M}_0)$, $(\mathcal{N}_0,\mathcal{M}_0)\prec (\mathcal{N}',\mathcal{M}')$, and $p$ is realized in $\mathcal{M}'=Z(U^{(\mathcal{N}',\mathcal{M}')})$.  
\end{proof}
\end{claim}

\begin{claim} If $\overline{b}\in\mathcal{N}_0$ and $p\in S_n(\overline{b})$, then there is $(\mathcal{N}',\mathcal{M}')\succ (\mathcal{N}_0,\mathcal{M}_0)$ with countable density character such that $p$ is realized in $\mathcal{N}'$.  
\begin{proof}
Consider the type $\{\phi(\overline{x},\overline{b})\dotminus \frac{1}{n}|\phi(\overline{x},\overline{b})\in p,n<\omega\}\cup Diag_{el}(\mathcal{N}_0,\mathcal{M}_0)$.  

As before, let $\Gamma(\overline{x})$ be a finite subset of this type and $\phi_1(\overline{x},\overline{b}),\ldots\phi_m(\overline{x},\overline{b})\in p$, $n_1,\ldots,n_m<\omega$ be such that $\Gamma(\overline{x})\subset \{\phi_1(\overline{x},\overline{b})\dotminus \frac{1}{n_1},\ldots,\phi_m(\overline{x},\overline{b})\dotminus \frac{1}{n_m}\}\cup Diag_{el}(\mathcal{N}_0,\mathcal{M}_0)$.  $p$ is realized in some $\mathcal{N}\succ \mathcal{N}_0$, so $\mathcal{N}\vDash\displaystyle{\inf_{\overline{x}}}\ \displaystyle{\max_{1\leq i\leq m}} \phi_i(\overline{x},\overline{b})=0$, which means $\mathcal{N}_0\vDash\displaystyle{\inf_{\overline{x}}}\ \displaystyle{\max_{1\leq i\leq m}} \phi_i(\overline{x},\overline{b})=0$.  Let $n=\displaystyle{\max_{1\leq i\leq m}} n_i$ and choose $\overline{x}\in\mathcal{N}_0$ such that $\displaystyle{\max_{1\leq i\leq m}} \phi_i(\overline{x},\overline{b}) \leq \frac{1}{n}$.  $\overline{x}$ realizes $\Gamma(\overline{x})$, so the type is finitely satisfiable in $\mathcal{N}_0$.  Thus we can choose $(\mathcal{N}',\mathcal{M}')$ with countable density character realizing the type.  $(\mathcal{N}',\mathcal{M}')\succ (\mathcal{N}_0,\mathcal{M}_0)$ and $p$ is realized in $\mathcal{N}'$.  
\end{proof}
\end{claim}

We build an elementary chain of models all with countable density character $(\mathcal{N}_0,\mathcal{M}_0)\prec (\mathcal{N}_1,\mathcal{M}_1)\prec\ldots$ such that 
\begin{enumerate}
\item[(i)] If $p\in S_n(T)$ is realized in $\mathcal{N}_{3i}$, then $p$ is realized in $\mathcal{M}_{3i+1}$.  
\item[(ii)] If $\overline{a},\overline{b},c\in\mathcal{M}_{3i+1}$ and $tp^{\mathcal{M}_{3i+1}}(\overline{a}) = tp^{\mathcal{M}_{3i+1}}(\overline{b})$, then there is $d\in\mathcal{M}_{3i+2}$ such that $tp^{\mathcal{M}_{3i+2}}(\overline{a},c) = tp^{\mathcal{M}_{3i+2}}(\overline{b},d)$.  
\item[(iii)] If $\overline{a},\overline{b},c\in\mathcal{N}_{3i+2}$ and $tp^{\mathcal{N}_{3i+2}}(\overline{a}) = tp^{\mathcal{N}_{3i+2}}(\overline{b})$, then there is $d\in\mathcal{N}_{3i+3}$ such that $tp^{\mathcal{N}_{3i+3}}(\overline{a},c) = tp^{\mathcal{N}_{3i+3}}(\overline{b},d)$. 
\end{enumerate}

(i) is possible by iterating the first claim.

For (ii), let $p= tp^{\mathcal{M}_{3i+1}}(c/\overline{a})$, and consider the type $\{\phi^U(\overline{b},x)\dotminus \frac{1}{n}|\phi(\overline{a},x)\in p,n<\omega\}\cup Diag_{el}(\mathcal{N}_{3i+1},\mathcal{M}_{3i+1})$.

As before, consider a finite subset of this type contained in some 

\noindent $\{\phi_1^U(\overline{b},x)\dotminus  \frac{1}{n_1},\ldots,\phi_m^U(\overline{b},x)\dotminus  \frac{1}{n_m}\}\cup Diag_{el}(\mathcal{N}_{3i+1},\mathcal{M}_{3i+1})$.  

Since $\phi_i(\overline{a},c) =0$ for $1\leq i\leq m$, $\mathcal{M}_{3i+1}\vDash\displaystyle{\inf_x} \displaystyle{\max_{1\leq i\leq m}} \phi_i(\overline{a},x) = 0$, so $\mathcal{M}_{3i+1}\vDash\displaystyle{\inf_x} \displaystyle{\max_{1\leq i\leq m}} \phi_i(\overline{b},x) = 0$ since  $tp^{\mathcal{M}_{3i+1}}(\overline{a}) = tp^{\mathcal{M}_{3i+1}}(\overline{b})$.  

Let $n = \displaystyle{\max_{1\leq i\leq m}} n_i$ and choose $x\in\mathcal{M}_{3i+1}$ such that $\displaystyle{\max_{1\leq i\leq m}} \phi_i(\overline{b},x)\leq \frac{1}{n}$.  Thus, this type is finitely satisfiable.  Choose $(\mathcal{N}_{3i+2},\mathcal{M}_{3i+2})\succ (\mathcal{N}_{3i+1},\mathcal{M}_{3i+1})$ with countable density character in which this type is realized.  Then $p$ is realized in $\mathcal{M}_{3i+2}$.  That is, there is $d\in\mathcal{M}_{3i+2}$ such that $tp^{\mathcal{M}_{3i+2}}(\overline{b},d) = tp^{\mathcal{M}_{3i+2}}(\overline{a},c)$.  

For (iii), we use the same argument as in (ii) to see that $\{\phi^U(\overline{b},x)\dotminus \frac{1}{n}|\phi(\overline{a},x)\in tp^{\mathcal{N}_{3i+2}}(c/\overline{a}),n<\omega\}\cup Diag_{el} (\mathcal{N}_{3i+2},\mathcal{M}_{3i+2})$ is finitely satisfiable in $\mathcal{N}_{3i+2}$ and apply the second claim.  

Let $(\mathcal{N},\mathcal{M}) = \displaystyle{\bigcup_{i<\omega}} (\mathcal{N}_i,\mathcal{M}_i)$.  Then $(\mathcal{N},\mathcal{M})$ has countable density character.  It is a Vaughtian pair, since (as we saw in Lemma \ref{countablevp}) being a Vaughtian pair is preserved by elementary extensions.  By (i), $\mathcal{M}$ and $\mathcal{N}$ realize the same types, and by (ii) and (iii), $\mathcal{M}$ and $\mathcal{N}$ are homogeneous.  Hence, by Proposition \ref{sametypes}, $\mathcal{M}\cong\mathcal{N}$. 
\end{proof}
\end{lemma}

\section{$(n,\epsilon,\delta)$-Quasi-Minimality}

In this section we propose a notion of Quasi-minimality in the continuous setting.  We find that, as is the case with other notions of minimality (\cite{noquez}), giving the correct characterization is a difficult task.  We wind up needing significant technical requirements in order to get our final result for the section, Lemma \ref{qmcor}, which we will need to prove Theorem \ref{kappamodel}.  

We will define $(n,\epsilon,\delta)$-quasi-minimality for $n<\omega$ and $0<\delta\leq \epsilon$ as a property of sets of the form $\{x\in\mathcal{M}|P(x)<K\}$ where $P$ is a definable predicate and $K>0$.  We say that a set of this form is $(n,\epsilon,\delta)$-quasi-minimal if the set has no countable $\epsilon$-net (is large with respect to $\epsilon$),  and for any subset of the form $\{x\in \mathcal{M}|Q(x)<r\}$ where $Q$ is a definable predicate and $r>0$, either this has a countable $\delta$-net (is small with respect to $\delta$), its complement is small with respect to $\delta$, or the set of elements in $\{x\in \mathcal{M}|P(x)<K\}$ which are $\frac{1}{n}$-away from every element of $\{x\in \mathcal{M}|Q(x)<r\}$ is small with respect to $\delta$.  In Proposition \ref{qmexistence}, we show that with an assumption of $\omega$-stability, for suitably chosen $\epsilon\geq \delta> 0$, we can find $P$ and $K$ such that $\{x\in\mathcal{M}|P(x)<K\}$ is $(n,\epsilon,\delta)$-quasi-minimal.  This proof uses a binary tree to contradict $\omega$-stability as in the classical case, but we see that it is necessary to specify that each set in our tree is at least $\frac{1}{n}$-apart, as we must guarantee that the uncountable set of types at the end does not have a countable dense subset with respect to the $d$-metric on types.  Using this, in Lemma \ref{qmlemma} we show that we can find a nested sequence of sets of this form which are quasi-minimal for strictly increasing values of $n$ such that their intersection has no countable $\epsilon$-net for some $\epsilon>0$.  This proof requires the third condition in the definition of $(n,\epsilon,\delta)$-quasi-minimality, which allows us to define uncountably many Cauchy sequences of elements of the nested sets such that at each stage of the intersection, we only remove countably many of them.  The limits of the uncountably many remaining Cauchy sequences will be in the intersection of these sets, and be $\epsilon$-apart from each other.  This intersection, which is quasi-minimal for all $n<\omega$ in some sense, will act similarly to a classical quasi-minimal set.  In Lemma \ref{qmcor} we use this intersection to define a type which is analogous to the collection of formulas which define large sets.

%
%

Let $T$ be a continuous $\omega$-stable theory.  Fix $\mathcal{M}\vDash T$.  For definable predicates $P$ over $\mathcal{M}$, we will just write $Z(P)$ to mean $Z(P^{\mathcal{M}})$ unless otherwise specified.  

In the classical setting, a definable set is \emph{quasi-minimal} if every definable subset is countable or co-countable. The following is an approximate version of quasi-minimality, which we give as a property of sets of the form $\{x\in \mathcal{M}|P(x)<K\}$ where $P$ is a definable predicate over $\mathcal{M}$ and $K>0$.  We will use $P<K$ as shorthand to denote the set $\{x\in\mathcal{M}|P(x)<K\}$.  

\begin{definition} For a definable predicate $P$ over $\mathcal{M}$ and $K>0$, $P<K$ is \emph{$(n,\epsilon,\delta)$-quasi-minimal in $\mathcal{M}$}, where $1\leq n<\omega$, $\epsilon>0$, and $\epsilon\geq \delta>0$, if $\{x\in\mathcal{M}|P(x)<K\}$ has no countable $\epsilon$-net and for every definable predicate $Q$ over $\mathcal{M}$, for every $r>0$, one of the following conditions holds:
\begin{enumerate}
\item[(1)] For every $r'$ such that $0<r'<r$, there exists $k'$ with $K\dotminus (r-r')<k'<K$ such that $\{x\in\mathcal{M}|Q(x)\leq r'\wedge P(x)\leq k'\}$ has a countable $\delta$-net.  

\item[(2)] $\{x\in\mathcal{M}|P(x)<K\wedge Q(x)\geq r\}$ has a countable $\delta$-net.
\item[(3)] There exists $r'$ such that $0<r'<r$ and for all $k'$ such that $K\dotminus (r-r')<k'<K$, $\{x\in\mathcal{M}|Q(x)< r'\wedge P(x)\leq k'\}$ has no countable $\delta$-net and $\{x\in\mathcal{M}|P(x)<K\wedge\forall y(Q(y)\leq r'\wedge P(y)\leq k'\rightarrow d(x,y)>\frac{1}{n})\}$ has a countable $\delta$-net.   

\end{enumerate}
\end{definition}

The idea is that  a set of the form $\{x\in \mathcal{M}|P(x)<K\}$ where $P$ is a definable predicate and $K>0$ is quasi-minimal if for any subset of the form $\{x\in \mathcal{M}|Q(x)<r\}$ where $Q$ is a definable predicate and $r>0$, either this set is small, its complement is small, or the set of elements in $\{x\in \mathcal{M}|P(x)<K\}$ which are $\frac{1}{n}$-away from every element of $\{x\in \mathcal{M}|Q(x)<r\}$ is small. 

Note that condition (1) holds if and only if $\{x\in \mathcal{M}|Q(x)<r\wedge P(x)<K\}$ has a countable $\delta$-net.  Also note that if $P<K$ is $(n,\epsilon,\delta)$-quasi-minimal, then for any $\epsilon'$ such that $\delta\leq \epsilon'\leq\epsilon$, $P<K$ is $(n,\epsilon',\delta)$-quasi-minimal, and for any $\delta'$ such that $\delta\leq \delta'\leq \epsilon$, $P<K$ is $(n,\epsilon,\delta')$-quasi-minimal.

We will now use $\omega$-stability to show the existence of an $(n,\epsilon,\delta)$-quasi-minimal subset of sets of the form $P<K$ for a fixed $n\in\omega$, and appropriately chosen $\epsilon$ and $\delta$.  Analogous to the classical setting, we will build a binary tree of sets of the form $P_\sigma <K_\sigma$ and use this to find uncountably many types which have no countable dense subset with respect to the $d$-metric on types (this is why we must specify $n$ in our characterization of quasi-minimality).  Note that because our sets are open, we will require some auxiliary $0<J_\sigma<K_\sigma$.

\begin{proposition}\label{qmexistence}  Let $T$ be $\omega$-stable and $\mathcal{M}\vDash T$.  Suppose $R$ is a definable predicate over $\mathcal{M}$ and $j>0$ is such that $\{x\in\mathcal{M}|R(x)<j\}$ has no countable $\epsilon$-net for some $\epsilon>0$.  Let $n<\omega$ such that $1\leq n$ and $\widehat{\epsilon}\in [0,\epsilon)$ be given.  Then there exists $\epsilon'$ such that $\widehat{\epsilon}<\epsilon'\leq\epsilon$, $P$ a definable predicate over $\mathcal{M}$, and $K>0$ such that $\{x\in\mathcal{M}|P(x)<K\}\subset \{x\in\mathcal{M}|R(x)<j\}$, and $P<K$ is $(n,\epsilon',\delta)$-quasi-minimal for all $\delta\in (\widehat{\epsilon}, \epsilon']$.  

\begin{proof}
Suppose not.  We will construct a binary tree of definable predicates to contradict $\omega$-stability.  

First note that there exists $j'\in (0,j)$ such that $Z(R\dotminus j')$ has no countable $\epsilon$-net.  If not, then $\{x\in\mathcal{M}|R(x)<j\}= \displaystyle{\bigcup_{\frac{1}{j}<n<\omega}} Z(R\dotminus (j\dotminus \frac{1}{n}))$ would have a countable $\epsilon$-net.  Choose $K_{\emptyset}$ such that $j'<K_{\emptyset}<j$.  Let $J_{\emptyset} = j$.  Let $\epsilon_{\emptyset}$ be such that $\widehat{\epsilon}<\epsilon_{\emptyset}\leq \epsilon$. 

Let $P_{\emptyset}=R$ and note that $\{x\in\mathcal{M}|P_{\emptyset}(x)<K_{\emptyset}\}$ contains $Z(R\dotminus j')$, so it has no countable $\epsilon_{\emptyset}$-net, and $\{x\in \mathcal{M}|P_{\emptyset}(x)<K_{\emptyset}\}\subset \{x\in \mathcal{M}|R(x)<j\}$.    

Suppose for some $\sigma\in 2^{<\omega}$ we have $P_{\sigma}$ a definable predicate over $\mathcal{M}$, and $0<K_{\sigma}<J_{\sigma}$ such that $\{x\in\mathcal{M}|P_{\sigma}(x)<K_{\sigma}\}$ has no countable $\epsilon_{\sigma}$-net for some $\epsilon_{\sigma}\in (\widehat{\epsilon},\epsilon]$, and  $\{x\in\mathcal{M}|P_\sigma(x)<J_{\sigma}\}\subset \{x\in R(x)<j\}$.  

We will find $P_{\sigma\widehat{\ }0}$ and $P_{\sigma\widehat{\ }1}$ definable predicates over $\mathcal{M}$, $\epsilon_{\sigma\widehat{\ }0}$ and $\epsilon_{\sigma\widehat{\ }1}$ in $(\widehat{\epsilon},\epsilon_{\sigma}]$, and $0<K_{\sigma\widehat{\ }0}<J_{\sigma\widehat{\ }0}$ and $0<K_{\sigma\widehat{\ }1}<J_{\sigma\widehat{\ }1}$ such that
\begin{itemize}
\item For all $x,y\in\mathcal{M}$, if $P_{\sigma\widehat{\ }0}(x)<J_{\sigma\widehat{\ }0}$ and $P_{\sigma\widehat{\ }1}(y)<J_{\sigma\widehat{\ }1}$, then $d(x,y)\geq\frac{1}{n}$.  
\end{itemize}
and for $i=0,1$,
\begin{itemize}
\item $\{x\in\mathcal{M}|P_{\sigma\widehat{\ }i}(x)<K_{\sigma\widehat{\ }i}\}$ has no countable $\epsilon_{\sigma\widehat{\ }i}$-net
\item For all $x\in\mathcal{M}$, $P_{\sigma\widehat{\ }i}(x)<J_{\sigma\widehat{\ }i}\Rightarrow P_{\sigma}(x)<K_{\sigma}$.  So $\{x\in\mathcal{M}|P_{\sigma\widehat{\ }i}<K_{\sigma\widehat{\ }i}\}\subset \{x\in\mathcal{M}|R(x)<j\}$.
\end{itemize}

First we will define $P_{\sigma\widehat{\ }1}$.  Since $\{x\in\mathcal{M}|P_\sigma(x)<K_\sigma\}\subset \{x\in\mathcal{M}|R(x)<j\}$ and $\widehat{\epsilon}<\epsilon_{\sigma}\leq \epsilon$, by assumption there exists $\delta\in (\widehat{\epsilon}, \epsilon_\sigma]$ such that $P_{\sigma}<K_{\sigma}$ is not $(n,\epsilon_{\sigma},\delta)$-quasi-minimal.  Since $\{x\in \mathcal{M}|P_{\sigma}(x)<K_{\sigma}\}$ has no countable $\epsilon_{\sigma}$-net, there is $\widehat{r}>0$ and $Q$ a definable predicate over $\mathcal{M}$ such that conditions (1), (2), and (3) from the definition of $(n,\epsilon_{\sigma},\delta)$-quasi-minimal all fail.  By the negation of (1), there exists $r'\in (0,\widehat{r})$ such that for all $k$ with $K_\sigma\dotminus (\widehat{r}-r')<k<K_\sigma$, 

\begin{equation}\label{one} \{x\in \mathcal{M}|Q(x)\leq r'\wedge P(x)\leq k\} \text{ has no countable $\delta$-net.}\tag{*}\end{equation}  

Then note that for all $r$ such that $r'<r<\widehat{r}$, for all $k$ such that $K_\sigma \dotminus (\widehat{r}-r')<k<K_\sigma$, $\{x\in \mathcal{M}|Q(x)< r\wedge P(x)\leq k\}$ has no countable $\delta$-net, since it contains 

\noindent $\{x\in \mathcal{M}|Q(x)\leq r'\wedge P(x)\leq k\}$.

By the negation of (3), for every $r$ such that $r'<r<\widehat{r}$, there exists $k_r$ with $K_\sigma\dotminus (\widehat{r}-r)<k_r<K_\sigma$ such that \begin{multline} \{x\in\mathcal{M}|P_\sigma(x)<K_\sigma\wedge \forall y(Q(y)\leq r\wedge P_\sigma(y)\leq k_r\rightarrow d(x,y)>\frac{1}{n})\}\\
 \text{ has no countable $\delta$-net}\tag{**}\end{multline} (since $\{x\in\mathcal{M}|Q(x)< r\wedge P_{\sigma}(x)\leq k_r\}$ has no countable $\delta$-net, because $K_{\sigma}\dotminus (\widehat{r}-r')<K_{\sigma}\dotminus (\widehat{r}-r)<k_r<K_{\sigma}$).  

Fix $r$ and $r_0$ such that $r'<r<r_0<\widehat{r}$, and let $k_0 = k_{r_0}$.  Observe that $$K_\sigma \dotminus (\widehat{r}-r')<K_{\sigma} \dotminus (\widehat{r}-r)<K_\sigma \dotminus (\widehat{r}-r_0)<k_0 < K_\sigma,$$ so we may choose $k'$ and $k$ such that $$K_\sigma \dotminus (\widehat{r}-r')<k'<K_{\sigma} \dotminus (\widehat{r}-r)<k<K_\sigma \dotminus (\widehat{r}-r_0)<k_0 < K_\sigma.$$  

We know from (*) that $\{x\in\mathcal{M}|Q(x)\leq r'\wedge P_\sigma(x)\leq k'\}$ has no countable $\delta$-net, and it is contained in $\{x\in\mathcal{M}|Q(x)<r\wedge P_\sigma(x)<k\}$, so by monotonicity this set has no countable $\delta$-net, and it is contained in $\{x\in\mathcal{M}|Q(x)\leq r_0\wedge P_\sigma(x)\leq k_0\}$.  We also know from (**) that $$\{x\in\mathcal{M}|P_{\sigma}(x)<K_\sigma\wedge \forall y(Q(y)\leq r_0\wedge P_{\sigma}(y)\leq k_0\rightarrow d(x,y)>\frac{1}{n})\}$$ has no countable $\delta$-net.  

Let $P_{\sigma\widehat{\ }1}(x) = \max(Q(x)\dotminus (r\dotminus k), P_{\sigma}(x)\dotminus (k\dotminus r))$ and $K_{\sigma\widehat{\ } 1} = \min(r,k)$.  So for all $x\in\mathcal{M}$, $P_{\sigma\widehat{\ }1}(x)<K_{\sigma\widehat{\ }1}$ if and only if $Q(x)<r$ and $P_\sigma(x)<k$.  

Thus, $\{x\in\mathcal{M}|P_{\sigma\widehat{\ }1}(x)<K_{\sigma\widehat{\ }1}\}$ has no countable $\delta$-net since it contains 

\noindent $\{x\in \mathcal{M}|Q(x)\leq r'\wedge P_{\sigma}(X)\leq k'\}$, and is contained in $\{x\in\mathcal{M}|Q(x)\leq r_0\wedge P_{\sigma}(x)\leq k_0\}$.  Let $\epsilon_{\sigma\widehat{\ }1} = \delta$, so $\widehat{\epsilon}<\epsilon_{\sigma\widehat{\ }1}\leq \epsilon$.  

Let $J_{\sigma\widehat{\ }1} = \min(r,k)+\min (r_0-r,k_0-k)$.  For $x\in \mathcal{M}$, if $P_{\sigma\widehat{\ }1}(x)<J_{\sigma\widehat{\ } 1}$,then $Q(x)<r_0$ and $P_\sigma(x)<k_0$.  

So since $k_0<K_\sigma$, if $P_{\sigma\widehat{\ }1}(x)<J_{\sigma\widehat{\ }1}$, then $P_{\sigma}(x)<K_\sigma$.  

Towards defining $P_{\sigma\widehat{\ }0}$, let $L(x) = \max(Q(x)\dotminus (r_0\dotminus k_0),P_\sigma(x)\dotminus (k_0\dotminus r_0))$ and $l = \min(r_0,k_0)>0$.  $L(x)\leq l$ if and only if $Q(x)\leq r_0$ and $P_{\sigma}(x)\leq k_0$.  

Rewritten in terms of $L$ and $l$,  $$\{x\in\mathcal{M}|P_\sigma(x)<K_\sigma \wedge \forall y (L(y)\leq l\rightarrow d(x,y)>\frac{1}{n})\}$$ has no countable $\delta$-net. 

Let $t = \min(K_\sigma,1-\frac{1}{n})>0$.  

Let $A(x) =P_\sigma(x)\dotminus (K_\sigma\dotminus (1-\frac{1}{n}))$.  Observe that $A(x)<t$ if and only if $P_\sigma(x)<K_\sigma$.

Let $B(x) =\displaystyle{\sup_y}\ \min (\frac{1-\frac{1}{n}}{1-l}(1-L(y)), 1-d(x,y))\dotminus ((1-\frac{1}{n})\dotminus K_\sigma)$ and observe that $B(x)<t$ if and only if for all $y$, either $\frac{1-\frac{1}{n}}{1-l}(1-L(y))<1-\frac{1}{n}$ or $1-d(x,y)<1-\frac{1}{n}$.  This is true if and only if for all $y$, either $L(y)>l$ or $d(x,y)>\frac{1}{n}$.  

That is, $B(x)<t$ if and only if for all $y$, $L(y)\leq l\rightarrow d(x,y)>\frac{1}{n}$.  (Observe that this complicated formula is necessary to express this implication since we cannot negate formulas in the continuous setting.)

Let $T(x) = \max(A(x), B(x))$.  Then $T(x)<t$ if and only if $x$ is in the set $\{x\in\mathcal{M}|P_\sigma(x)<K_\sigma \wedge \forall y (L(y)\leq l\rightarrow d(x,y)>\frac{1}{n})\}$, so $\{x\in\mathcal{M}|T(x)<t\}$ has no countable $\delta$-net.  

Then there exists $t'\in (0,t)$ such that $\{x\in\mathcal{M}|T(x)\leq t'\}$ has no countable $\delta$-net.  Let $P_{\sigma\widehat{\ }0} = T$ and let $K_{\sigma\widehat{\ }0}$ be such that $t'<K_{\sigma\widehat{\ }0}<t$.  Let $J_{\sigma\widehat{\ }0} = t$.  $\{x\in\mathcal{M}|P_{\sigma\widehat{\ }0}(x)<K_{\sigma\widehat{\ }0}\}$ has no countable $\delta$-net,  Let $\epsilon_{\sigma\widehat{\ }0}=\delta$, so $\widehat{\epsilon}<\epsilon_{\sigma\widehat{\ }0}\leq \epsilon$.  If $P_{\sigma\widehat{\ }0}(x)<J_{\sigma\widehat{\ }0}$ then $P_{\sigma}(x)<K_{\sigma}$ for all $x\in\mathcal{M}$.  

Finally let $x\in\mathcal{M}$ be such that $P_{\sigma\widehat{\ }0}(x)<J_{\sigma\widehat{\ }0}$ and $y\in\mathcal{M}$ be such that $P_{\sigma\widehat{\ }1}(y)<J_{\sigma\widehat{\ }1}$.  Then $T(x)<t$, and $Q(y)< r_0$ and $P_{\sigma}(y)< k_0$.  Since $T(x)<t$, for all $y\in\mathcal{M}$, if $Q(y)\leq r_0$ and $P_{\sigma}(y)\leq k_0$, then $d(x,y)>\frac{1}{n}$.  Thus, $d(x,y)>\frac{1}{n}$. 

Let $A$ be the union over $\sigma\in 2^{<\omega}$ of the (countable) set of parameters from $\mathcal{M}$ used to define each $P_{\sigma}$.  This is a countable union of countable sets, so $A$ is countable.  

Let $\tau\in 2^{\omega}$ be given.  For $n<\omega$, let $m^{\tau}_n = \frac{K_{\tau|_n} + J_{\tau|_n}}{2}$.  So $K_{\tau|_n}<m^{\tau}_n<J_{\tau|_n}$.  Let $p_{\tau} = \{P_{\tau|_n}(x)\dotminus m^{\tau}_n:n<\omega\}$. We will show that $p_{\tau}$ is a consistent type.  Note that for $x\in \mathcal{M}$, if $P_{\tau|_{n+1}}(x)<J_{\tau|_{n+1}}$, then $P_{\tau|_n}(x)<K_{\tau|_n}$.  Consider a finite subset of $p_{\tau}$ contained in $\{P_{\tau|_0}(x)\dotminus m^{\tau}_0,\ldots, P_{\tau|_n}(x)\dotminus m^{\tau}_n\}$.  Then since we know that $\{x\in \mathcal{M}|P_{\tau|_n}(x)\leq m^\tau_n\}\supset \{x\in\mathcal{M}|P_{\tau|_n}(x)<K_{\tau|_n}\}$ has no countable $\epsilon_{\tau|_n}$-net, it is non-empty, so this finite part is realized.  Thus, $p_{\tau}$ is a consistent partial type over $A$.  

Let $\tau\neq \tau'$ in $2^\omega$ be given and let $\sigma\in 2^{<\omega}$ with $|\sigma|=n$ be such that $\tau|_n=\sigma=\tau'|_n$ and (without loss of generality) $\tau|_{n+1} = \sigma\widehat{\ }0$ and $\tau'|_{n+1} = \sigma\widehat{\ }1$.

Let $\mathcal{N}\succ\mathcal{M}$ be an elementary extension of $\mathcal{M}$ with elements $a_\tau\vDash p_\tau$ and $a_{\tau'}\vDash p_{\tau'}$.  

We know that for all $x\in\mathcal{M}$, if $P_{\sigma\widehat{\ }0}(x)<J_{\sigma\widehat{\ }0}$ and $P_{\sigma\widehat{\ }1}(y)<J_{\sigma\widehat{\ }1}$, then $d(x,y)>\frac{1}{n}$.  Let $\theta(x,y) = \max(P_{\sigma\widehat{\ }0}(x)\dotminus (J_{\sigma\widehat{\ }0}\dotminus J_{\sigma\widehat{\ }1}), P_{\sigma\widehat{\ }1}(y)\dotminus (J_{\sigma\widehat{\ }1}\dotminus J_{\sigma\widehat{\ }0}))$ and let $m =\min(J_{\sigma\widehat{\ }0},J_{\sigma\widehat{\ }1})$.  So $\theta(x,y)<m$ if and only if $P_{\sigma\widehat{\ }0}(x)<J_{\sigma\widehat{\ }0}$ and $P_{\sigma\widehat{\ }1}(y)<J_{\sigma\widehat{\ }1}$.

So for all $x,y\in\mathcal{M}$, $\theta(x,y)<m$ implies $d(x,y)>\frac{1}{n}$, so $d(x,y)\geq \frac{1}{n}$.  Equivalently, for all $x,y\in \mathcal{M}$, $\theta(x,y)\geq m$ or $d(x,y)\geq \frac{1}{n}$.  We can express this with $$\mathcal{M}\vDash \displaystyle{\sup_x}\ \displaystyle{\sup_y}\ \min(m\dotminus \theta(x,y),\frac{1}{n}\dotminus d(x,y))=0,$$ so by elementarity, 

$$\mathcal{N}\vDash \displaystyle{\sup_x}\ \displaystyle{\sup_y}\ \min(m\dotminus \theta(x,y),\frac{1}{n}\dotminus d(x,y))=0.$$ 

Let $x\vDash tp(a_{\tau}/A)$ and $y\vDash tp(a_{\tau'}/A)$.  Then $P_{\sigma\widehat{\ }0}(x) = P_{\sigma\widehat{\ }0}(a_\tau) \leq m^\tau_{n+1} < J_{\sigma\widehat{\ }0}$ and $P_{\sigma\widehat{\ }1}(y)= P_{\sigma\widehat{\ }1}(a_{\tau'})\leq m^{\tau'}_{n+1} <J_{\sigma\widehat{\ }1}$.  So $\theta(x,y)<m$, which means $m\dotminus \theta(x,y)>0$.  Thus, we must have $\frac{1}{n}\dotminus d(x,y)=0$, so $d(x,y)\geq \frac{1}{n}$.  Since $x\vDash tp(a_{\tau}/A)$ and $y\vDash tp(a_{\tau'}/A)$ were arbitrary, the distance between these types is at least $\frac{1}{n}$.  

Thus, $\{tp(a_{\tau}/A)|\tau\in2^{\omega}\}$ is an uncountable collection of complete types over the countable set $A$ whose sets of realizations are pairwise $\frac{1}{n}$-apart, contradicting $\omega$-stability (since it cannot have a countable dense subset with respect to the d-metric on types).
\end{proof}
\end{proposition}

Next we will show that it is possible to find a set which is, in some sense, quasi-minimal for all $n\in\omega$.  We do this by using Lemma \ref{qmexistence} to find a nested sequence of sets which are quasi-minimal for increasingly large $n$.  Using the third condition in the definition of quasi-minimality, we find uncountably many Cauchy sequences which live in the intersection of all of these sets whose limits have no countable dense subset.  This intersection is analogous to a classical quasi-minimal set.

For the remainder of this section, fix $T$ an $\omega$-stable continuous theory, $\mathcal{M}\vDash T$ and $\epsilon>0$ such that $\mathcal{M}$ has no countable $\epsilon$-net.  

For $n<\omega$, let $m_n = \lceil \frac{1}{\epsilon} \rceil 2^{n+3}$.  Then $\displaystyle{\sum_{n<\omega}} \frac{2}{m_n} \leq \displaystyle{\sum_{n<\omega}} \frac{2}{(\frac{1}{\epsilon})2^{n+3}} = \frac{\epsilon}{2} \displaystyle{\sum_{n<\omega}} \frac{1}{2^{n+1}} = \frac{\epsilon}{2}$.  

Let $\epsilon_n = \epsilon - \displaystyle{\sum_{k<n}} \frac{2}{m_k}$ and let $\epsilon' = \epsilon - \displaystyle{\sum_{n<\omega}} \frac{2}{m_n}$.  

So $\epsilon=\epsilon_0>\epsilon_1>\ldots>\epsilon'>0$ and $\epsilon_n\rightarrow\epsilon'$.  

The definable predicates in the next lemma will be used to define the type in Lemma \ref{qmcor}, which is used in turn to prove Lemma \ref{technicallemma}.

\begin{lemma}\label{qmlemma} There is a strictly increasing sequence of integers $(i_n:n<\omega)$ with $i_n\geq n$ such that we can choose a sequence of definable predicates $P_n$ over $\mathcal{M}$ and $K_n>0$ such that $P_n<K_n$ is $(m_{n},\epsilon_{i_n}, \delta)$-quasi-minimal for all $\delta\in (\epsilon', \epsilon_{i_n}]$, $\{x\in \mathcal{M}|P_n(x)<K_n\}\supset \{x\in \mathcal{M}|P_{n+1}(x)<K_{n+1}\}$, and if $B_n\subset \{x\in \mathcal{M}|P_n(x)<K_n\}$ are open sets such that $B_n$ has a countable $\epsilon_{i_{n+1}}$-net, then $\displaystyle{\bigcap_{n<\omega}} (\{x\in \mathcal{M}|P_n(x)<K_n\}\setminus B_n)$ has no countable $\frac{\epsilon'}{2}$-net.  

\begin{proof}

$\mathcal{M} = \{x\in\mathcal{M}|d(x,x)<1\}$, which has no countable $\epsilon$-net.  By Proposition \ref{qmexistence}, there is $P_0$ a definable predicate over $\mathcal{M}$, $\widehat{\epsilon}\in (\epsilon',\epsilon)$, and $K_0>0$ such that $P_0<K_0$ is $(m_0,\widehat{\epsilon}, \delta)$-quasi-minimal for all $\delta\in (\epsilon',\widehat{\epsilon})$.  Let $i_0\geq 0$ be such that $\epsilon'<\epsilon_{i_0}<\widehat{\epsilon}$.  Then $P_0<K_0$ is $(m_0,\epsilon_{i_0},\delta)$-quasi-minimal for all $\delta\in (\epsilon', \epsilon_{i_0}]$.  

Then there exists $K'_0\in (0,K_0)$ such that $\{x\in\mathcal{M}|P_0(x)<K'_0\}$ has no countable $\epsilon_{i_0}$-net.

Suppose for some $n<\omega$ we have a definable predicate $P_n$ over $\mathcal{M}$, $K_n>0$, and $i_n\geq n$ such that $P_n<K_n$ is $(m_n,\epsilon_{i_n},\delta)$-quasi-minimal for all $\delta\in (\epsilon',\epsilon_{i_n}]$ and $K'_n\in (0,K_n)$ such that $\{x\in \mathcal{M}|P_n(x)<K'_n\}$ has no countable $\epsilon_{i_n}$-net. 

By Proposition \ref{qmexistence}, there is a definable predicate $P_{n+1}$ over $\mathcal{M}$ and $K_{n+1}>0$ such that for some $\widehat{\epsilon}\in (\epsilon',\epsilon_{i_n}]$, $P_{n+1}<K_{n+1}$ is $(m_{n+1},\widehat{\epsilon},\delta)$-quasi-minimal for all $\delta(\epsilon', \widehat{\epsilon}]$.  Choose $i_{n+1}\geq \max(n+1,i_n)$ such that $\epsilon_{i_{n+1}}<\widehat{\epsilon}$.  So $P_{n+1}<K_{n+1}$ is $(m_{n+1},\epsilon_{i_{n+1}},\delta)$-quasi-minimal for all $\delta\in (\epsilon',\epsilon_{i_{n+1}}]$ and $\{x\in\mathcal{M}|P_{n+1}(x)<K_{n+1}\}\subset \{x\in\mathcal{M}|P_n(x)<K'_n\}\subset\{x\in\mathcal{M}|P_n(x)<K_n\}$.  

Choose $K'_{n+1}<K_{n+1}$ such that $\{x\in\mathcal{M}|P_{n+1}(x)<K'_{n+1}\}$ has no countable $\epsilon_{i_{n+1}}$-net.

So since $P_n<K_n$ is $(m_n,\epsilon_{i_n},\epsilon_{i_{n+1}})$-quasi-minimal, and condition (1) does not hold for $P_{n+1}$ and $K'_{n+1}$, either condition (2) holds or condition (3) holds.

If condition (2) holds, then $\{x\in\mathcal{M}|P_{n+1}(x)\geq K'_{n+1}\wedge P_n(x)<K_n\}$ has a countable $\epsilon_{i_{n+1}}$-net, so $\{x\in \mathcal{M}|P_{n+1}(x)\geq K'_{n+1}\wedge P_n(x)<K'_n\}$ has a countable $\epsilon_{i_{n+1}}$-net.  Let $$A_n = \{x\in\mathcal{M}|P_{n+1}(x)\geq K'_{n+1}\wedge P_n(x)<K'_n\}$$ and say $A_n$ is of Case I. 

If not, then condition (3) must hold.  Let $r\in (0,K'_{n+1})$ be such that for all $k$ with $K_n\dotminus (K'_{n+1}-r)<k<K_n$, $\{x\in \mathcal{M}|P_{n+1}(x)<r\wedge P_n(x)\leq k\}$ has no countable $\epsilon_{i_{n+1}}$-net and $\{x\in\mathcal{M}|P_n(x)<K_n\wedge \forall y(P_{n+1}(y)\leq r \wedge P_n(y) \leq k \rightarrow d(x,y)>\frac{1}{m_n}\}$ has a countable $\epsilon_{i_{n+1}}$-net.  

Choose $k$ such that $\max(K'_n,K_n\dotminus (r-r'))<k<K_n$.  

So $\{x\in\mathcal{M}|P_n(x)<K_n\wedge \forall y(P_{n+1}(y)\leq r \wedge P_n(y) \leq k \rightarrow d(x,y)>\frac{1}{m_n})\} = \{x\in\mathcal{M}|P_n(x)<K_n\wedge \forall y(P_{n+1}(y)\leq r\rightarrow d(x,y)>\frac{1}{m_n})\}$, since $P_{n+1}(y)\leq r\Rightarrow P_{n+1}(y)<K_{n+1}\Rightarrow P_n(y)\leq K'_n<k$.  Let $A_n = \{x\in \mathcal{M}|P_n(x)<K'_n \wedge \forall y(P_{n+1}(y)\leq r\rightarrow d(x,y)>\frac{1}{m_n})\}\subset \{x\in\mathcal{M}|P_n(x)<K_n\wedge \forall y(P_{n+1}(y)\leq r\rightarrow d(x,y)>\frac{1}{m_n})\}$ which has a countable $\epsilon_{i_{n+1}}$-net, so $A_n$ has a countable $\epsilon_{i_{n+1}}$-net.  We say that $A_n$ is of Case II. 

If $x\in \{x\in\mathcal{M}|P_n(x)<K'_n\}\setminus A_n$, then there exists $y\in \{x\in \mathcal{M}|P_{n+1}(y)\leq r\}\subset \{x\in\mathcal{M}|P_{n+1}(y)<K'_{n+1}\}$ with $d(x,y)\leq \frac{1}{m_n}$.  

So for every $n<\omega$ we have $P_n$ a definable predicate over $\mathcal{M}$, $0<K'_n<K_n$ such that $P_n<K_n$ is $(m_{n},\epsilon_{i_n},\delta)$-quasi-minimal for all $\delta\in (\epsilon', \epsilon_{i_n}]$, $\{x\in\mathcal{M}|P_n(x)<K'_n\}$ has no countable $\epsilon_{i_{n}}$-net, and $P_{n+1}(x)<K_{n+1}$ implies $P_n(x)<K'_n$ for all $x\in \mathcal{M}$.  

We also have $A_n\subset \{x\in \mathcal{M}|P_n(x)<K_n\}$ such that $A_n$ has a countable $\epsilon_{i_{n+1}}$-net.  

Now let $B_n$ be given as in the statement of the lemma.  Let $C_n = A_n\cup (B_n\cap \{x\in \mathcal{M}|P_n(x)<K'_n\})$, so $C_n$ has a countable $\epsilon_{i_{n+1}}$-net.

Now let $(a^0_n:i\in I_0)$ be a collection of pairwise $\epsilon_0$-apart elements of $\{x\in\mathcal{M}|P_0(x)<K'_0\}$ with $I_0$ an uncountable index set.  

Suppose we have defined $L_0,\ldots,L_{n-1}\subset I_0$ all countable, and $I_n = I_0\setminus (L_0\cup\ldots\cup L_{n-1})$, and we have $(a^n_i:i\in I_n)$ which are pairwise $\epsilon_n$-apart in $\{x\in \mathcal{M}|P_n(x)<K'_n\}$.  

$C_n$ has a countable $\epsilon_{i_{n+1}}$-net, so since $i_{n+1}\geq n+1>n$, it has a countable $\epsilon_n$-net.  So $\{a^n_i:i\in I_n\}\cap C_n$ must be countable, since otherwise $C_n$ would have no countable $\epsilon_n$-net since $\epsilon_n>\epsilon_{i_{n+1}}$.  

Let $L_{n}\subset I_n$ be the set of indices $i$ such that such that $a^n_i\in C_n$.  $L_n$ is countable, let $I_{n+1} = I_n\setminus L_n$.  

Let $i\in I_{n+1}$ be given.  $a^n_i\notin B_n$ and  $a^n_i\notin A_n$.  If $A_n$ is of Case I, then $a^n_i$ is in 

\noindent $\{x\in\mathcal{M}|P_{n+1}(x)<K'_{n+1}\}$, so let $a^{n+1}_i = a^n_i$.  

If $A_n$ is of Case II, there exists $y\in\mathcal{M}$ such that $P_{n+1}(y)<K'_{n+1}$ and $d(a^n_i,y)\leq \frac{1}{m_n}$.  Let $a^{n+1}_i=y$, so $P_{n+1}(a^{n+1}_i)<K'_{n+1}$.  

In both cases, $a^{n+1}_i\in \{x\in \mathcal{M}|P_{n+1}(x)<K'_{n+1}\}$, and $d(a^{n+1}_i,a^n_i)\leq \frac{1}{m_n}$.  Since for  $i,j\in I_{n+1}$, $d(a^n_i,a^n_j)>\epsilon_n$, $d(a^{n+1}_i,a^{n+1}_j)>\epsilon_n - \frac{2}{m_n} = \epsilon_{n+1}$.

Let $I = I_0\setminus \displaystyle{\bigcup_{n<\omega}}L_n$, so $I$ is still uncountable since every $L_n$ is countable. 

For each $i\in I$, the sequence $a^n_i$ is Cauchy, so  since $\mathcal{M}$ is complete, it converges to some $a_i$.  Then, for $n<\omega$, for all $m>n$, $a^m_i \in \{x\in \mathcal{M}|P_m(x)<K'_m\}\setminus B_m \subset \{x\in \mathcal{M}|P_m(x)\leq K'_m\}\setminus B_m$, so since this is a closed set, their limit $a_i\in \{x\in\mathcal{M}|P_m(x)\leq K'_m\}\setminus B_m$.  This is true of all $m<\omega$, so $a_i\in \displaystyle{\bigcap_{n<\omega}}( \{x\in\mathcal{M}|P_n(x)\leq K'_n\}\setminus B_n) \subset\displaystyle{\bigcap_{n<\omega}} (\{x\in\mathcal{M}|P_n(x)<K_n\}\setminus B_n)$.  

Finally, for $i\neq j$, since $d(a^n_i,a^n_j)> \epsilon_n\geq\epsilon'$ for all $n$, $d(a_i,a_j)\geq\epsilon'>\frac{\epsilon'}{2}$.  Hence, 

\noindent $\displaystyle{\bigcap_{n<\omega}} (\{x\in\mathcal{M}|P_n(x)<K_n\}\setminus B_n)$ has no countable $\frac{\epsilon'}{2}$-net, as required.  

\end{proof}
\end{lemma}

Next we will use the $P_n<K_n$ from the previous lemma to define a type $p$ such that $\phi=0$ is in $p$ if and only if for every $r>0$, $\{x\in\mathcal{M}|\phi(x)>r\}$ has a countable $\epsilon$-net for an appropriately chosen $\epsilon>0$.  That is, the type consists of formulas whose zero sets are approximately ``co-countable''.  

\begin{lemma}\label{qmcor} Let $P_n<K_n$ be as in the previous lemma, so $P_n<K_n$ is $(m_n,\epsilon_{i_n},\epsilon_m)$-quasi-minimal for all $m\geq i_n$.  

Let $p$ be a collection of $\mathcal{L}_{\mathcal{M}}$-conditions such that $\phi=0$ is in $p$ if and only if for every $r>0$ there exists $N_r<\omega$ such that for all $n>N_r$, $\{x\in \mathcal{M}|P_n(x)<K_n\wedge \phi(x)\geq r\}$ has a countable $\epsilon_{i_{n+1}}$-net.   

Then $p$ is a consistent complete type over $\mathcal{M}$. 

\begin{proof}
First we will show that $p$ is complete by showing that for every $\mathcal{L}_{\mathcal{M}}$-formula $\phi$, there is $r\in [0,1]$ such that $|\phi-r|=0$ is in $p$.

Let $\phi$ an $\mathcal{L}_{\mathcal{M}}$-formula be given.  Suppose $\phi=0$ is not in $p$.  Then there exists $r>0$ such that for all $N<\omega$ there exists $n>N$ such that  $\{x\in \mathcal{M}|P_n(x)<K_n\wedge \phi(x)\geq r\}$ has no countable $\epsilon_{i_{n+1}}$-net.  

Let $r$ be the supremum of the set of $r>0$ witnessing this.  Let $r'>r$ be given.  Then $r'$ does not witness this, so there exists $N<\omega$ such that for all $n>N$, $\{x\in\mathcal{M}|\phi(x)\geq r'\wedge P_n(x)<K_n\}$ has a countable $\epsilon_{i_{n+1}}$-net.  

Suppose $|\phi-r|=0$ is not in $p$.  Then there exists $s>0$ such that for all $N<\omega$ there is $n>N$ such that $\{x\in\mathcal{M}|P_n(x)<K_n\wedge|\phi(x)-r|\geq s\}$, which is equal to $$\{x\in \mathcal{M}|P_n(x)<K_n\wedge\phi(x)\leq r\dotminus s\}\cup \{x\in\mathcal{M}|\phi(x)\geq r\dot{+}s\wedge P_n(x)<K_n\},$$ has no countable $\epsilon_{i_{n+1}}$-net.  

So by choosing $N$ sufficiently large, we know that for all $m>N$, there exists $n>m$ such that $\{x\in \mathcal{M}|P_n(x)<K_n\wedge\phi(x)\leq r\dotminus s\}$ has no countable $\epsilon_{i_{n+1}}$-net, since $\{x\in\mathcal{M}|\phi(x)\geq r\dot{+}s\wedge P_n(x)<K_n\}$ does have a countable $\epsilon_{i_{n+1}}$-net, because it is contained in $\{x\in \mathcal{M}|P_n(x)<K_n\wedge \phi(x)\geq r\}$.

Let $k_0$ and $k_1$ be such that $r\dotminus s<k_0<k_1<r$.   Let $N$ be such that $\frac{1}{m_N}<\Delta_\phi(k_1-k_0)$, where $\Delta_{\phi}$ is the modulus of uniform continuity for $\phi$. 

Choose $n_1>N$ so that $\{x\in \mathcal{M}|P_n(x)<K_n\wedge \phi(x)>k_1\}$  has no countable $\epsilon_{i_{n_1+1}}$ net (since it contains $\{x\in \mathcal{M}|P_{n_1}(x)<K_{n_1}\wedge \phi(x)\geq r\}$), and choose $n_2>n_1$ sufficiently large so that $\{x\in \mathcal{M}|P_{n_2}(x)<K_{n_2}\wedge \phi(x)<k_0\}$ has no countable $\epsilon_{i_{n_2+1}}$-net (since it contains $\{x\in \mathcal{M}|P_{n_2}(x)<K_{n_2}\wedge \phi(x)\leq r\dotminus s\}$).  Then $\{x\in \mathcal{M}|P_{n_1}(x)<K_{n_1}\wedge \phi(x)<k_0\}$ has no countable $\epsilon_{i_{n_2}+1}$-net, since by Lemma \ref{qmlemma}, $P_{n_2}(x)<K_{n_2}\Rightarrow P_{n_1}(x)<K_{n_1}$, and $\{x\in \mathcal{M}|P_{n_1}(x)<K_{n_1}\wedge \phi(x)>k_1\}$ has no countable $\epsilon_{i_{n_2+1}}$-net since $\epsilon_{i_{n_2+1}}<\epsilon_{i_{n_1+1}}$.  Let $\delta = \epsilon_{i_{n_2+1}}$ and $n=n_1$.  Note that $\delta >\epsilon'$.

So we have $\{x\in \mathcal{M}|P_n(x)<K_n\wedge \phi(x)<k_0\}$ which has no countable $\delta$-net as well as $\{x\in \mathcal{M}|P_n(x)<K_n\wedge \phi(x)>k_1\}$ with no countable $\delta$-net.  $P_n<K_n$ is $(m_n,\epsilon_{i_n},\delta)$-quasi-minimal, so since conditions (1) and (2) fail for $\phi$ and $k_0$, there is $r'$ such that $0<r'<k_0$ and for all $k'$ such that $K_n-(k_0-r')<k'<K_n$, $\{x\in \mathcal{M}|\phi(x)<r'\wedge P_n(x)\leq k'\}$ has no countable $\delta$-net (so in particular, it is not empty) and $$\{x\in \mathcal{M}|P_n(x)<K_n\wedge \forall y(P_n(y)\leq k'\wedge \phi(y)\leq r'\rightarrow d(x,y)>\frac{1}{m_n})\}$$ has a countable $\delta$-net.  

Thus, there must be some $x\in \{x\in \mathcal{M}|P_n(x)<K_n\wedge \phi(x)>k_1\}$ such that $x$ is not in $\{x\in \mathcal{M}|P_n(x)<K_n\wedge \forall y(P_n(y)\leq k'\wedge \phi(y)\leq r'\rightarrow d(x,y)>\frac{1}{m_n})\}$ since the latter set has a countable $\delta$-net, and the former does not.  

Then there exists $y\in \{x\in \mathcal{M}|\phi(y)\leq r'\wedge P_n(x)\leq k'\} \subset \{x\in \mathcal{M}|\phi(x)<k_0 \wedge P_n(x)<K_n\}$ such that $d(x,y)\leq \frac{1}{m_n}\leq\frac{1}{m_N}$.

Since $d(x,y)\leq \frac{1}{m_N}<\Delta_\phi(k_1-k_0)$, $|\phi(x)-\phi(y)|\leq k_1-k_0$, but this is a contradiction, since $\phi(y)<k_0$ and $\phi(x)>k_1$, so $|\phi(x)-\phi(y)|>k_1-k_0$. 

So $|\phi-r|=0$ must be in $p$. 

Thus, $p$ is a complete type.

To see that $p$ is consistent, let $\phi_1=0,\ldots,\phi_m=0$ in $p$ be given.  Let $r>0$ be given.  Choose $N$ sufficiently large such that for all $n>N$ and $i$ such that $1\leq i\leq m$, the set $\{x\in\mathcal{M}|\phi_i(x)\geq r\wedge P_n(x)<K_n\}$ has a countable $\epsilon_{i_{n+1}}$-net.  Then the union of these sets has a countable $\epsilon_{i_{n+1}}$-net, so $\{x\in \mathcal{M}|\max(\phi_1(x),\ldots,\phi_n(x))\geq r\wedge P_n(x)<K_n\}$ has a countable $\epsilon_{i_{n+1}}$-net.  

Hence $\max(\phi_1,\ldots,\phi_n)=0$ is in $p$.

Finally, for $\phi=0$ in $p$, $\{\phi(x)\leq \frac{1}{n}:n<\omega\}$ is finitely satisfiable, since for every $n<\omega$, there is $m<\omega$ such that $\{x\in \mathcal{M}|P_m(x)<K_m(x)\wedge \phi(x)\geq \frac{1}{n}\}$ has a countable $\epsilon_{i_{n+1}}$-net, so $\{x\in\mathcal{M}|P_m(x)<K_m(x)\wedge \phi(x)<\frac{1}{n}\}$ has no countable $\epsilon_{i_{n+1}}$-net, and in particular, is not empty.  Thus, this partial type is finitely satisfiable, so it is consistent.  

Hence, since $p$ is closed under finite conjunctions, by compactness, $p$ is consistent.  
\end{proof}
\end{lemma}

\section{Vaught's Two-Cardinal Theorem}

In this section we will prove Vaught's Two-Cardinal Theorem, as well as a partial converse of the theorem which requires an additional assumption of $\omega$-stability.  Ultimately, these will be necessary for our result about uncountable categoricity.

For the first result (Theorem \ref{vaughtstwocardinaltheorem}) we prove it in a similar fashion to the classical analogue.  However, for the partial converse (Theorem \ref{kappamodel}), we need Lemma \ref{qmcor}, which required our development of approximate quasi-minimality in the continuous setting. While the argument via quasi-minimality follows the proof in the classical case, as we saw in the previous section, the development of a notion of quasi-minimality in this setting required a substantial departure from the classical setting.  

\begin{definition} Let $\kappa>\lambda\geq \aleph_0$.  We say that an $\mathcal{L}$-theory $T$ has a $(\kappa,\lambda)$-model if there is $\mathcal{M}\vDash T$ and $P$ a definable distance predicate in $\mathcal{M}$ such that $\mathcal{M}$ has density character $\kappa$ and $Z(P^{\mathcal{M}})$ has density character $\lambda$.  \end{definition}

Note that if $\mathcal{M}\vDash T$ is a $(\kappa,\lambda)$-model, then $T$ is not $\kappa$-categorical (by a straight forward compactness argument, we can show that there exists a model of $T$ with density character $\kappa$ such that every non-compact zero set of a definable predicate has density character $\kappa$).

We will begin by proving the following theorem of Vaught in the continuous setting: If $T$ has a $(\kappa,\lambda)$-model for some $\kappa>\lambda\geq \aleph_0$, then $T$ has an $(\aleph_1,\aleph_0)$-model.

\begin{lemma} If $T$ has a $(\kappa,\lambda)$-model, then $T$ has a Vaughtian pair.
\begin{proof} Let $\mathcal{N}\vDash T$ be a $(\kappa,\lambda)$-model with $P$ a definable distance predicate over some countable $A\subset\mathcal{N}$ such that $Z(P^{\mathcal{N}})$ has density character $\lambda$.  Then by Downward L\"{o}wenheim Skolem (Proposition 7.3 in \cite{mtfms}), there is $\mathcal{M}\prec \mathcal{N}$ containing $Z(P^{\mathcal{N}})$ and $A$ with density character $\lambda$.  So $\mathcal{M}\neq\mathcal{N}$, and $Z(P^{\mathcal{M}})=Z(P^{\mathcal{N}})$.  Thus, $(\mathcal{N},\mathcal{M})$ is a Vaughtian pair of models of $T$.
\end{proof}
\end{lemma}

\begin{proposition}\label{aleph1} If $T$ has a Vaughtian pair, then $T$ has an $(\aleph_1,\aleph_0)$-model.
\begin{proof}
By Lemma \ref{countablevp} and Lemma \ref{homogvp}, there is $(\mathcal{N},\mathcal{M})$ a Vaughtian pair with countable density character such that $\mathcal{M}$ and $\mathcal{N}$ are homogeneous models of $T$ realizing the same types, and $\mathcal{M}\cong\mathcal{N}$.  

Let $P$ be a definable distance predicate in $\mathcal{M}$ such that $Z(P^{\mathcal{M}}) = Z(P^{\mathcal{N}})$ and $Z(P^{\mathcal{M}})$ is not compact.  

We build an elementary chain $(\mathcal{N}_\alpha:\alpha<\omega_1)$ such that $\mathcal{N}_\alpha \cong\mathcal{N}$ and $(\mathcal{N}_{\alpha+1},\mathcal{N}_\alpha) \cong (\mathcal{N},\mathcal{M})$ as $\mathcal{L}^*$-structures, so in particular, $Z(P^{\mathcal{N}_{\alpha+1}}) = Z(P^{\mathcal{N}_\alpha})$ and $\mathcal{N}_{\alpha}\neq\mathcal{N}_{\alpha+1}$.  

Let $\mathcal{N}_0 = \mathcal{N}$, and for $\alpha$ a limit ordinal, let $\mathcal{N}_\alpha$ be the completion of $\displaystyle{\bigcup_{\beta<\alpha}} \mathcal{N}_\beta$.  Since $\mathcal{N}_\alpha$ is the union of homogeneous models isomorphic to $\mathcal{N}$, $\mathcal{N}_\alpha$ is homogeneous and realizes the same types as $\mathcal{N}$, so by Proposition \ref{sametypes}, $\mathcal{N}_\alpha \cong \mathcal{N}$.  

Finally, given $\mathcal{N}_\alpha$, $\mathcal{N}_\alpha\cong\mathcal{N}\cong\mathcal{M}$ so there is an elementary extension 

\noindent $\mathcal{N}_{\alpha+1}\succneqq\mathcal{N}_\alpha$ such that $(\mathcal{N}_{\alpha+1},\mathcal{N}_\alpha)\cong(\mathcal{N},\mathcal{M})$.  So $\mathcal{N}_{\alpha+1}\cong\mathcal{N}$.  

Let $\mathcal{N}^*$ be the completion of $\displaystyle{\bigcup_{\alpha<\omega_1}} \mathcal{N}_\alpha$.  There is $x\in\mathcal{N}\setminus\mathcal{M}$ such that $d(x,\mathcal{M}) =\delta>0$.  So $\mathcal{N}_{\alpha+1}\vDash \displaystyle{\inf_x}\displaystyle{\inf_{y\in Z(U)}} |d(x,y)-\delta| = 0$ since $U$ is a definable distance predicate. This is realized in $\mathcal{N}_\alpha$, so it is realized in  $\mathcal{N}_{\alpha+1}$.  So there are $\omega_1$ $x$'s which are $\delta$ apart in $\mathcal{N}^*$.  Thus, $\mathcal{N}^*$ has density character $\aleph_1$.   $|\mathcal{N}^*| = \aleph_1$ since for each $\alpha<\omega_1$, $|\mathcal{N}_{\alpha} |\leq \aleph_1$.

If $\mathcal{N}^*\vDash P(\overline{a})$, then $\overline{a}\in \mathcal{M}$, which has countable density character, so $Z(P^{\mathcal{N}^*})$ has countable density character.

Hence, $\mathcal{N}^*$ is an $(\aleph_1,\aleph_0)$-model of $T$.  

\end{proof}
\end{proposition}

Thus, we have proved Vaught's Two-Cardinal theorem for continuous logic:

\begin{theorem}\label{vaughtstwocardinaltheorem} If $T$ has a $(\kappa,\lambda)$-model where $\kappa>\lambda\geq \aleph_0$, then $T$ has an $(\aleph_1,\aleph_0)$-model.  \end{theorem}

Now we will show a partial converse of this theorem with the additional assumption that $T$ is $\omega$-stable, namely that if there is an $(\aleph_1,\aleph_0)$-model of $T$, then for any $\kappa\geq\aleph_1$, there is a $(\kappa,\aleph_0)$-model of $T$.  This will require the results about quasi-minimal sets from the previous section.  

Let $T$ be a continuous $\omega$-stable theory.

The next lemma tells us that in an $\omega$-stable theory $T$, for $\mathcal{M}\vDash T$ and a countable type which satisfies certain technical requirements, if this type is realized in an elementary extension of $\mathcal{M}$, then it has an approximate realization in $\mathcal{M}$.  

\begin{lemma}\label{technicallemma}
Suppose $T$ is $\omega$-stable, $\mathcal{M}\vDash T$, and $\mathcal{M}$ is not separable.  There is a proper elementary extension $\mathcal{N}$ of $\mathcal{M}$ with the following property:  suppose $\Gamma(\overline{w})$ is a countable type over $\mathcal{M}$ such that for every $\gamma(\overline{w})\in \Gamma(\overline{w})$, $\Delta_\gamma(\epsilon)=\epsilon$ for all $\epsilon\in[0,1]$, where $\Delta_\gamma$ is the modulus of uniform continuity for $\gamma$.  Then if $\Gamma(\overline{w})$ is realized in $\mathcal{N}$, for every $\epsilon>0$ there is $\overline{b}_\epsilon\in\mathcal{M}$ such that $\gamma(\overline{b}_\epsilon)\leq \epsilon$ for all $\gamma(\overline{w})\in\Gamma(\overline{w})$.  
\begin{proof}
Let $p$ be the type described in Lemma \ref{qmcor}.  

Let $\mathcal{M}'\succ \mathcal{M}$ be an elementary extension of $\mathcal{M}$ containing $c\notin\mathcal{M}$ realizing $p$.  Since $p$ is a complete type over $\mathcal{M}$, $p = tp(c/\mathcal{M})$.  

Since $T$ is $\omega$-stable, there is a prime model $\mathcal{N}\prec\mathcal{M}'$ containing $\mathcal{M}\cup\{c\}$ such that every type over $\mathcal{M}\cup \{c\}$ realized in $\mathcal{N}$ is principal.  So for $q$ which is realized in $\mathcal{N}$, the set of realizations of $q$ in $\mathcal{N}$ is a definable set. 

Let $\Gamma(\overline{w})$ be a countable type over $\mathcal{M}$ such that for all $\gamma\in \Gamma$, $\Delta_\gamma(\epsilon)=\epsilon$ for $\epsilon\in[0,1]$.  Let $\overline{b}\in\mathcal{N}$ realize $\Gamma$.  Then there is a predicate $Q(\overline{w},c)$ which is definable over $\mathcal{M}$ such that $Q(\overline{w},c)$ is the distance predicate for the set of realizations of $\Gamma(\overline{w})$.  In particular, if $\gamma(\overline{w})\in \Gamma(\overline{w})$, and $\overline{b}'$ is such that $Q(\overline{b}',c)= \epsilon$, then for all $\delta>0$ there is $\overline{b}'_\delta\in \mathcal{N}$ realizing $\Gamma$ with $d(\overline{b}'_{\delta},\overline{b}')\leq \epsilon+\delta$.  Thus, since $\Delta_\gamma(\epsilon+\delta) = \epsilon+\delta$, and since $\gamma(\overline{b}'_\delta)=0$, $\gamma(\overline{b}')=|\gamma(\overline{b}')-\gamma(\overline{b}'_\delta)|\leq \epsilon+\delta$.  Since this is true for all $\delta>0$, $\gamma(\overline{b}')\leq \epsilon$.  Thus, $\gamma(\overline{b}')\leq Q(\overline{b}',c)$.  

So, since $Q(\overline{b},c)=0$ and $c\vDash p$, $\displaystyle{\inf_{\overline{w}}}\ Q(\overline{w},x)$ is in $p$ and $\displaystyle{\sup_{\overline{w}}}\ \gamma(\overline{w})\dotminus Q(\overline{w},x)$ is in $p$ for every $\gamma\in \Gamma$.  

Let $\Theta$ be the type $\{\displaystyle{\inf_{\overline{w}}}\ Q(\overline{w},x)\}\cup\{\displaystyle{\sup_{\overline{w}}}\ \gamma(\overline{w})\dotminus Q(\overline{w},x):\gamma\in\Gamma\}$.  So $\Theta\subset p$ is countable.  Let $\theta_0(x),\theta_1(x),\ldots$ enumerate $\Theta$.  

Let $j<\omega$ and $k<\omega$ such that $k\geq 1$ be given.  $\theta_j=0$ is in $p$, so let $n<\omega$ be the least such that $\{x\in\mathcal{M}|P_n(x)<K_n\wedge \theta_j(x)\geq \frac{1}{k}\}$ has a countable $\epsilon_{i_{n+1}}$-net.

Let $B^{(j,k)}_n = \{x\in\mathcal{M}|P_n(x)<K_n\wedge \theta_j(x)> \frac{1}{k}\}$, which also has a countable $\epsilon_{i_{n+1}}$-net and is open.

For $n'\neq n$, let $B^{(j,k)}_{n'} = \emptyset$.  So we have defined $B^{(j,k)}_n$ for all $j$, $k\geq 1$, and $n$. 

For $n<\omega$, let $B_n = \displaystyle{\bigcup_{j<\omega,1\leq k<\omega}}\ B^{(j,k)}_n$. So $B_n$ has a countable $\epsilon_{i_{n+1}}$-net, and is a union of open sets, so it is open.   

Let $A = \displaystyle{\bigcap_{n<\omega}}(\{x\in\mathcal{M}|P_n(x)<K_n\}\setminus B_n)$, where $P_n<K_n$ are as in Lemma \ref{qmlemma}.  By Lemma \ref{qmlemma}, $A$ has no countable $\frac{\epsilon'}{2}$-net.  So in particular, $A\neq \emptyset$.

Let $c'\in A$.  Then $\theta_j(c')\leq\frac{1}{k}$ for all $j<\omega$ and $k<\omega$ with $k\geq 1$, so $\theta_j(c')=0$ for all $j<\omega$.  Thus, $c'$ realizes $\Theta$.  

So $\displaystyle{\inf_{\overline{w}}}\ Q(\overline{w},c') = 0$ and $\displaystyle{\sup_{\overline{w}}}\ \gamma(\overline{w})\dotminus Q(\overline{w},c') = 0$ for all $\gamma\in \Gamma$. 

Let $\epsilon>0$ be given.  Then we can choose $\overline{b}_\epsilon\in\mathcal{M}$ such that $Q(\overline{b}_\epsilon,c')\leq \epsilon$.  Hence, for $\gamma\in \Gamma$, $\gamma(\overline{b}_\epsilon) \leq \epsilon$, as required.

\end{proof}
\end{lemma}

\begin{theorem}\label{kappamodel} Suppose $T$ is $\omega$-stable and that there is an $(\aleph_1,\aleph_0)$-model of $T$.  If $\kappa\geq \aleph_1$, then there is a $(\kappa,\aleph_0)$-model of $T$.  
\begin{proof}
Let $\mathcal{M}\vDash T$ with $|\mathcal{M}|\geq\aleph_1$ and uncountable density character, and let $P$ be a definable distance predicate such that $Z(P^{\mathcal{M}})$ has a countable dense subset $\{m_i:i<\omega\}$.  

Let $Q(v)$ be the (definable) predicate $dist(v,Z(P))=\displaystyle{\inf_{y\in Z(P)}} d(v,y)$.  Observe that $Z(Q)=Z(P)$ and the modulus of uniform continuity for $Q$ is $\Delta_Q(\epsilon)=\epsilon$ for $\epsilon\in[0,1]$.  For $n<\omega$, let $\Gamma_n(v)$ be the type $\{Q(v)\}\cup\{\frac{1}{n}\dotminus  d(v,m_i):i<\omega\}$.  Observe that $\Gamma_n(v)$ is countable and $\Delta_\gamma(\epsilon)=\epsilon$ for all $\epsilon\in[0,1]$ and $\gamma\in \Gamma_n$, so $\Gamma_n$ satisfies the conditions of the previous lemma.  Let $\mathcal{N}\succneqq\mathcal{M}$ be as in the lemma. 

Suppose for some $n<\omega$, $\Gamma_n(v)$ is realized in $\mathcal{N}$.  Let $\epsilon$ be such that $0<\epsilon<\frac{1}{4n}$.  By the lemma, there is $b\in\mathcal{M}$ such that $Q(b)\leq \epsilon$ and $\frac{1}{n}\dotminus d(b,m_i)\leq \epsilon$ for all $i<\omega$.  Since $Q(b)\leq \epsilon$, we can choose $y\in Z(P^{\mathcal{M}})$ such that $d(b,y)\leq 2\epsilon$.  Since $(m_i:i<\omega)$ is dense in $Z(P^{\mathcal{M}})$, we can choose $m_i$ such that $d(y,m_i)\leq \epsilon$.  Then $d(b,m_i)\leq d(b,y)+d(y,m_i)\leq 2\epsilon+\epsilon = 3\epsilon<\frac{3}{4n}$, but $d(b,m_i)\geq \frac{1}{n}-\epsilon>\frac{3}{4n}$. $\Rightarrow\Leftarrow$

So for every $v\in\mathcal{N}$, if $P(v)=0$, then for every $n<\omega$ there is $m_i$ with $d(v,m_i)\leq \frac{1}{n}$, or else $v$ would realize $\Gamma_n$.  Hence, $Z(P^{\mathcal{N}})$ has the same countable dense subset as $Z(P^{\mathcal{M}})$.

Note that there is $x\in\mathcal{N}\setminus\mathcal{M}$ with $dist(x,\mathcal{M}) =: \delta>0$.  

Iterating this construction, we build an elementary chain $(\mathcal{M}_\alpha:\alpha<\kappa)$ such that $\mathcal{M}_0=\mathcal{M}$, $\mathcal{M}_{\alpha+1}\neq\mathcal{M}_\alpha$ and $Z(P^{\mathcal{M}_{\alpha+1}})$ has the same countable dense subset as $Z(P^{\mathcal{M}_0})$ and there exists $x\in \mathcal{M}_{\alpha+1}\setminus\mathcal{M}_{\alpha}$ with $dist(x,\mathcal{M}_{\alpha})=\delta>0$.  For $\alpha$ a limit, let $\mathcal{M}_{\alpha}$ be the completion of $\displaystyle{\bigcup_{\beta<\alpha}}\mathcal{M}_{\beta}$.  

If $\mathcal{N}$ is the completion of $\displaystyle{\bigcup_{\alpha<\kappa}} \mathcal{M}_\alpha$, $\mathcal{N}$ has density character $\kappa$, since it has $\kappa$-many elements which are $\delta$-apart, so $\mathcal{N}$ is a $(\kappa,\aleph_0)$-model of $T$.
\end{proof}
\end{theorem}

\section{Uncountable Categoricity}

In classical logic, the Baldwin-Lachlan characterization of uncountable categoricity says that a theory $T$ is uncountably categorical if and only if $T$ is $\omega$-stable and has no Vaughtian pairs.  Though it has been shown that $\omega$-stability and the absence of Vaughtian pairs are not sufficient conditions for uncountable categoricity (\cite{hanson2020strongly}), here we use the results of the previous section to prove the forward direction of this theorem in the continuous setting.

\begin{definition} For an infinite cardinal $\kappa$, a continuous theory $T$ is \emph{$\kappa$-categorical} if for any $\mathcal{M},\mathcal{N}\vDash T$ with density character $\kappa$, $\mathcal{M}\cong\mathcal{N}$.  \end{definition}

\begin{theorem}\label{forwarddirection} If a continuous theory $T$ is $\kappa$-categorical for some $\kappa\geq \aleph_1$, then $T$ is $\omega$-stable and has no Vaughtian pairs.  
\begin{proof}
Let $\kappa\geq \aleph_1$ and suppose $T$ is $\kappa$-categorical.  By Theorem 5.2 in \cite{benyaacov}, $T$ is $\omega$-stable. Suppose $T$ has a Vaughtian pair.  Then by Proposition \ref{aleph1}, $T$ has an $(\aleph_1,\aleph_0)$-model.  So by Theorem \ref{kappamodel}, $T$ has a $(\kappa,\aleph_0)$-model.  Using a compactness argument, we can show that $T$ has a model with density character $\kappa$ such that every non-compact zero set of a definable predicate has density character $\kappa$.  Thus, $T$ cannot be $\kappa$-categorical. 
\end{proof}
\end{theorem}

\section{Examples}
In this section we give an example of a Vaughtian pair of models of the theory of the Urysohn Sphere, and an example of a Vaughtian pair of models of the randomization of a theory.  By Theorem \ref{forwarddirection}, this tells us that these theories cannot be $\kappa$-categorical for any $\kappa>\aleph_0$.  

\subsection{Urysohn Sphere}

Let $\mathfrak{U}$ denote the Urysohn sphere, which is the unique (up to isomorphism) universal complete separable metric space of diameter $1$ in the ``empty'' language (which only has a symbol $d$ for the metric).  Let $\Theta_n$ be the collection of formulas of the form $\displaystyle{\max_{1\leq i<j\leq n}}\ |d(x_i,x_j)-r_{i,j}|$ where $r_{i,j}>0$ are such that this is a possible distance configuration of $n$ distinct points in a metric space of diameter $1$ (that it, does not violate the triangle inequality).  For $\theta\in \Theta_{n+1}$, if $\theta=\displaystyle{\max_{1\leq i<j\leq n+1}}|d(x_i,x_j)-r_{i,j}|$ where $r_{i,j}>0$, let $\theta|_n = \displaystyle{\max_{1\leq i<j\leq n}}|d(x_i,x_j)-r_{i,j}|$.  

Clearly, for every $\theta\in \Theta_{n+1}$, for every $\epsilon>0$ there exists $\delta = \delta(\epsilon)>0$ such that for all $a_1,\ldots,a_n\in \mathfrak{U}$, if $\theta|_n(a_1,\ldots,a_n)<\delta$, then there exists $a_{n+1}\in\mathfrak{U}$ such that $\theta(a_1,\ldots,a_n,a_{n+1})\leq \epsilon$.  

Let $T_{\mathfrak{U}}$ be the collection of conditions of the form $$\displaystyle{\sup_{x_1}}\ldots\displaystyle{\sup_{x_n}}\ \displaystyle{\inf_{x_{n+1}}} \min(\frac{\epsilon}{1-\delta}(1-\theta|_n(x_1,\ldots,x_n)),\theta(x_1,\ldots,x_n,x_{n+1}))\dotminus \epsilon $$ for $\theta\in \Theta_{n+1}$.  

This is just another way of saying that for all $x_1,\ldots,x_n\in\mathfrak{U}$, there exists $x_{n+1}\in\mathfrak{U}$ such that if $\theta|_n(x_1,\ldots,x_n)<\delta$, then $\theta(x_1,\ldots,x_n,x_{n+1})\leq \epsilon$.  




Fact 5.1 in \cite{urysohn} tells us that $\mathfrak{U}$ is the only separable model of $T_{\mathfrak{U}}$, so $T_{\mathfrak{U}}$ is $\aleph_0$-categorical, and thus, complete.   By Proposition 5.3 in \cite{urysohn}, $T_{\mathfrak{U}}$ admits quantifier elimination, so it is model complete.  

Let $y\in\mathfrak{U}$ and let $\mathcal{M} = \mathfrak{U}\setminus\{x\in\mathfrak{U}|d(x,y)<\frac{1}{8}\}$.  $\mathcal{M}$ is a complete separable metric space with diameter 1.  We will show that $\mathcal{M}\vDash T_{\mathfrak{U}}$.  

Let $\theta\in \Theta_{n+1}$ be given.  Suppose $t_1,\ldots,t_n\in [\frac{1}{8},1]$ are such that 

\noindent $\max(\theta|_n(x_1,\ldots,x_n), \displaystyle{\max_{1\leq i\leq n}}|d(x_i,y)-t_i|)$ is a possible configuration of $n+1$ points.  In other words, $\theta|_n(x_1,\ldots,x_n)$ is a possible configuration of $n$ points in $\mathcal{M}$.  

Let $r_i = d(x_i,x_{n+1})$ in $\theta$ for $1\leq i\leq n$.  It is easy to check that by the triangle inequality, for all $1\leq i,j\leq n$, $|t_i-r_i|\leq t_j+r_j$. 

Let $t_{n+1} = \displaystyle{\min_{1\leq j\leq n}} (t_j + r_j)$.  Observe that $t_{n+1}>\frac{1}{8}$ since $r_j>0$ and $t_j\geq \frac{1}{8}$, and let $s = t_{n+1} - \frac{1}{8} >0$.  

Note that $t_{n+1} \geq \displaystyle{\max_{1\leq i\leq n}} |t_i-r_i|$.  

It is easy to verify that $\max(\theta(x_1,\ldots,x_{n+1}),\displaystyle{\max_{1\leq i\leq n+1}}\ |d(x_i,y)-t_i|) =0$ does not violate the triangle inequality, which means that this is in $\Theta_{n+1}$.

So, given $\epsilon>0$, there exists $\delta>0$ such that $T_{\mathfrak{U}}$ says that for all $x_1,\ldots, x_n,y$ there exists $x_{n+1}$ such that if $\theta|_n(x_1,\ldots,x_n)<\delta$ and $|d(x_i,y)-t_i|<\delta$ for $1\leq i\leq n$, then $\theta(x_1,\ldots,x_n,x_{n+1})\leq \min(\epsilon,s)$ and $|d(x_{n+1},y)-t_{n+1}|\leq \min(\epsilon,s)$.  

Let $x_1,\ldots,x_n\in\mathcal{M}$ be given and assume $\theta|_n(x_1,\ldots,x_n)<\delta$.  Then, there exists $x_{n+1}$ such that $\theta(x_1,\ldots,x_n,x_{n+1})\leq \min(\epsilon,s)\leq \epsilon$, and $|d(x_{n+1},y)-t_{n+1}|\leq \min(\epsilon,s)\leq s$.  So $d(x_{n+1},y)\geq t_{n+1}-s= \frac{1}{8}$.  Hence, $x_{n+1}\in\mathcal{M}$.

So $\mathcal{M}\vDash \displaystyle{\sup_{x_1}}\ldots\displaystyle{\sup_{x_n}}\ \displaystyle{\inf_{x_{n+1}}} \min(\frac{\epsilon}{1-\delta}(1-\theta|_n(x_1,\ldots,x_n)),\theta(x_1,\ldots,x_n,x_{n+1}))\dotminus \epsilon$.  Thus, $\mathcal{M}\vDash T_{\mathfrak{U}}$.  

Since $\mathcal{M}\subset \mathfrak{U}$, by model completeness, $\mathcal{M}\prec\mathfrak{U}$.  Since $y\notin \mathcal{M}$, $\mathcal{M}\neq\mathfrak{U}$.  

Choose $x\in \mathfrak{U}$ with $d(x,y) = \frac{1}{2}$ (this exists since $\mathfrak{U}\vDash T_{\mathfrak{U}}$).  Let $P(v) = d(x,v)\dotminus \frac{1}{8}$.  Observe that $Z(P)$ is a definable set, since $dist(v,Z(P)) = d(v,y)\dotminus \frac{1}{8}$.  

So if $v\in \mathfrak{U}$ and $v\in Z(P)$, then $d(x,v)\leq \frac{1}{8}$, so $d(v,y)\geq d(x,y) - d(x,v) \geq \frac{1}{2}-\frac{1}{8} = \frac{3}{8} >\frac{1}{8}$.  Thus, $v\in\mathcal{M}$.  So $Z(P^{\mathfrak{U}}) = Z(P^{\mathcal{M}})$

Further note that $Z(P^{\mathfrak{U}})$ is not compact, since compact subsets of $\mathfrak{U}$ have empty interior (Corollary 2.13 in \cite{isaac}).  

Hence, $(\mathfrak{U},\mathcal{M})$ is a Vaughtian pair of models of $T_{\mathfrak{U}}$.

\subsection{Randomizations}

In this section we will apply the result of the previous section to show that for a (classical or continuous) theory $T$, its randomized theory $T^R$ is not $\kappa$-categorical for any $\kappa\geq \aleph_1$.

Randomizations of theories are, for the most part, model theoretically similar to the original theory.  In particular, in \cite{BK}, Ben Yaacov and Keisler showed that $\omega$-categoricity, $\omega$-stability, and stability are preserved, in \cite{BenNIP}, Ben Yaacov shows that NIP is preserved, and in \cite{AK}, Andrews and Keisler show that $T$ has a prime model if and only if $T^R$ has a prime model, and that if the original theory $T$ is $\aleph_1$-categorical, then $T^R$ has at most countably many separable models.  Not all model theoretic properties are preserved, in \cite{bensimple}, Ben Yaacov shows that the randomization of a simple, stable structure is not simple.  And we see here that since we can always find a Vaughtian pair of models of $T^R$, by Theorem \ref{forwarddirection}, the randomization of an uncountably categorical theory is not uncountably categorical.

A randomization of a model $\mathcal{M}$ of a (classical or continuous) theory $T$ is a two-sorted continuous structure with a sort $\mathcal{K}$ whose elements are random elements of $\mathcal{M}$, and a sort $\mathcal{B}$ whose elements are events in an underlying probability space.  We assume familiarity with Keisler randomizations viewed as metric structures, but we will recall the basics here for the reader's convenience.  For a more complete introduction, see Section 2 in \cite{AK}.  

Fix a classical or continuous countable language $L$ and a complete $L$-theory $T$.  The randomization signature $L^R$ is the continuous language with sorts $\mathcal{K}$ and $\mathcal{B}$, an $n$-ary function symbol $\llbracket\phi(\cdot)\rrbracket:\mathcal{K}^n\rightarrow \mathcal{B}$ for each $L$-formula $\phi$, and the Boolean operations $\top,\perp,\sqcup,\sqcap,\neg$ of sort $\mathcal{B}$.  If $(\Omega,\mathcal{B},\mu)$ is our underlying probability space, for $\overline{f}$ of sort $\mathcal{K}$, $\llbracket \phi(\overline{f})\rrbracket = \{\omega\in\Omega|\mathcal{M}\vDash \phi(\overline{f}(\omega)\}$, so its measure can be though of as the probability that $\overline{f}$ is in $\phi(\mathcal{M})$.

For $f,g$ of sort $\mathcal{K}$, $d(f,g)$ is the measure of $\llbracket f\neq g\rrbracket$ and for $A,B$ of sort $\mathcal{B}$, $d(A,B)$ is the measure of their symmetric difference. 

Here, we will restrict our attention to the case when our underlying atomless finitely additive probability algebra is $([0,1),\mathcal{L},\lambda)$ where $\mathcal{L}$ are the Borel subsets of $[0,1)$ and $\lambda$ is Lebesgue measure.  Our $L^R$ pre-structures will be of the form $(\mathcal{M}^{[0,1)},\mathcal{L})$ where $\mathcal{M}\vDash T$.  By Fact 2.5 in \cite{AK}, for any $\mathcal{M}\vDash T$, if we identify elements which are distance $0$ from each other, $(\mathcal{M}^{[0,1)},\mathcal{L})$ is an $L^R$-structure.  

\begin{fact}[Theorem 2.1 in \cite{BK}] There is a unique complete $L^R$-theory $T^R$, called the randomized theory of $T$, such that for every model $\mathcal{M}$ of $T$, $(\mathcal{M}^{[0,1)},\mathcal{L})$ is a pre-model of $T$.  \end{fact}

We will find a Vaughtian pair of models of this $T^R$.

Let $\mathcal{M},\mathcal{N}\vDash T$ be such that $\mathcal{M}\precneqq \mathcal{N}$.    Let $\mathcal{K}_0,\mathcal{K}_1\vDash T^R$ be the structures obtained from $(\mathcal{M}^{[0,1)}, \mathcal{L})$ and  $(\mathcal{N}^{[0,1)}, \mathcal{L})$.

For $x\in \mathcal{N}$, let $f_x$ in $\mathcal{K}_1$ denoted the constant function $f(t) = x$ for all $t\in [0,1)$.  

By Remark 2.3 in \cite{AK},  $\mathcal{K}_0\prec \mathcal{K}_1$.  Further note that $\mathcal{K}_0\neq \mathcal{K}_1$, since, for example, for $x\in \mathcal{N}\setminus\mathcal{M}$, $f_x$ is in the random variable sort of $\mathcal{K}_1$, but not in the random variable sort of $\mathcal{K}_0$.  

Fix $a,b\in\mathcal{M}$ with $d(a,b)>0$ (or just $a\neq b$ if $\mathcal{M}$ is first order). Let $\phi(g)$ be $|\lambda(\llbracket g=f_a\rrbracket) + \lambda (\llbracket g=f_b\rrbracket) -1|$.  So for $g\in \mathcal{N}^{[0,1)}$, $\phi(g)=0$ if and only if the range of $g$ is in $\{a,b\}$.  So if $g\in Z(\phi^{\mathcal{K}_1})$, then $g\in Z(\phi^{\mathcal{K}_0})$, since $a,b\in\mathcal{M}$.  That is, $\phi$ has the same zero set in $\mathcal{K}_0$ and $\mathcal{K}_1$.  To see that $(\mathcal{K}_1,\mathcal{K}_0)$ is a Vaughtian pair, we need to show that $Z(\phi)$ is definable and not compact.

Let $\epsilon>0$ be given.  If $\phi(g)<\epsilon$, then $\lambda (\{t\in [0,1)|g(t)\notin\{a,b\}\}) <\epsilon$.  Define $$g'(t) = \begin{cases} 
      a & g(t)=a \\
      b & g(t)=b \\
      a & otherwise 
   \end{cases}$$

$g'\in Z(\phi)$, and $d(g,g') = \lambda (\{t\in [0,1)|g(t)\neq g'(t)\}) \leq \epsilon$.  Thus, $dist(g,Z(\phi))\leq \epsilon$.  So by Proposition \ref{defdistpred}, $Z(\phi)$ is definable.

Finally, let $0<\epsilon<\frac{1}{2}$ be given.  Suppose $f_1,\ldots,f_m\in Z(\phi)$ are centers of a finite $\epsilon$-net of $Z(\phi)$.  Let $A_i = \{t\in [0,1)|f_i(t) = a\}\in \mathcal{L}$.  For $\sigma\in 2^m$, let $A_\sigma = \displaystyle{\bigcap_{1\leq i\leq m}} A_i^{\sigma(i)}$ where $A_i^0 = [0,1)\setminus A_i$ and $A_i^1 = A_i$.  So $\{A_\sigma|\sigma\in 2^m\}$ is a partition of $[0,1)$ into elements of $\mathcal{L}$ such that each $f_i$ is constant on each $A_\sigma$.  

By atomlessness, for each $\sigma\in 2^{m}$, there is $B_\sigma\in \mathcal{L}$, $B_\sigma\subset A_\sigma$ such that $\lambda(B_\sigma) = \frac{1}{2}\lambda(A_\sigma)$.  Define \[g(t) = \begin{cases} a & t\in B_\sigma$ for some $\sigma\in 2^m\\
b & t\in A_\sigma\setminus B_\sigma$ for some $\sigma \in 2^m \end{cases}\]

$\phi(g) =0$, since $\lambda(\llbracket g=f_a\rrbracket) + \lambda(\llbracket g=f_b\rrbracket) = \displaystyle{\sum_{\sigma\in 2^m}} \mu(B_\sigma) + \displaystyle{\sum_{\sigma\in 2^m}} \mu(A_\sigma\setminus B_\sigma) = 1$.  And for any $1\leq i\leq m$, $d(g,f_i) = \displaystyle{\sum_{\sigma \in 2^m}} \frac{1}{2}\mu(A_\sigma) = \frac{1}{2}>\epsilon$.  So there is no finite $\epsilon$-net of $Z(\phi)$.  

Thus, $(\mathcal{K}_1,\mathcal{K}_0)$ is a Vaughtian pair of models of $T^R$.  Hence, by Theorem \ref{forwarddirection}, $T^R$ is not $\kappa$-categorical for any uncountable $\kappa$.

\bibliography{twocardinal}{}
\bibliographystyle{amsalpha}
\end{document}